# GEORG CANTOR FROM ST. PETERSBURG.
# CHILDHOOD AND HISTORY OF THE FAMILY.
# ARCHIVAL RESEARCH

***G.I. Sinkevich, St.-Petersburg State University of architecture and civil engineering, St. Petersburg, Russia.***
*galina.sinkevich@gmail.com*

Georg Cantor (1845-1918) was born, and spent the first 11 years of his life in St. Petersburg. The present lecture is devoted to his childhood and his family. Most of these documents were not available before and are now published for the first time. (Picture 1. Generation table - at the end of the article).

We know the history of his family from his letters, his German archives. There are excellent works by A. Frenkel, H. Meschkovski, I. Grattan-Guinness, I.-W. Dauben, W. Purkert and H.J. Ilgauds, Ann-Marie Decallot about Cantor, his life and works. But some information about his origin and family was unknown to Georg Cantor himself. In Petersburg archived I have find many interesting information about his family` life.

His family was remarkable. There were eight violinists, two professors (violin and law), ladies of the Imperial Court, merchants and painters.

Georg Ferdinand Ludwig Philipp Cantor was born in St. Petersburg on March, 3 (February, 19 in the old style), 1845, in the family of Georg Woldemar Cantor and Maria Cantor-Boehm. The weather was cold – minus 8–10 degrees Reaumur (10–13 degrees Celsius), there was a gentle breeze, and the sky was covered by clouds.

Georg was the first long-awaited child, who was born in the third year after the marriage of his parents. His father wrote the following in 1851: "The elder one, seven-year-old, was named Georg Ferdinand Ludwig Philipp, and he was gifted by the nature, with an aspiration for order, which prevailed upon any other aspirations, and he was sanguine person in nature" [1]. (Picture 2. Fragment of the letter).

His father, Georg Woldemar Cantor, merchant and broker, moved to St. Petersburg from Copenhagen. His mother, Maria Cantor-Boehm, was born in St. Petersburg in a family of violinists.

The family lived in the Vasilyevskiy Island, in the 11th line, in the house of merchant Tranchel[1] [2]. Now there is a memorial plaque in the yard. (Pictures 3, 4,5).

There was a beautiful garden near the house, the garden created by the gardener Gerdes, and the children of the Cantor's family used to play in the garden.

From this house, the father of the family, Georg Woldemar Cantor, used to go to the Strelka (the eastern tip of Vasilyevsky Island) where he was a broker at the Stock Exchange, and via the Blagovestchensky bridge to the Isaac Square, to the Trading House of Sarepta Society[2], where he was a collaborator. (Picture 6. Stock Exchange). In the department of manuscripts of the Russian National Library (RNL), a letter written by Georg Woldemar Cantor in 1851 is still available, in which letter he described "severe bustle of his life", exhausting competition, his apathy, his children, guests in his house, cold winter evenings, boiled beef for supper, and boiling samovar".

Georg was baptized in St. Catherine Lutheran Church not far from his house. Apropose, Leonard Euler was a parishioner of this Church in XVIII century. (Picture 7. St. Catherine Lutheran Church).

In 1846, the second boy was born, named Ludwig Gustav (Lelya).

---

[1] On October 10, 2011, a memorial plaque was placed in the yard of the house, at the current address of St. Petersburg, 11th line of Vasilyevsky Island, No. 24.

[2] The current address is St. Petersburg, Yakubovich Street, No. 24.

From the letter written by Georg Woldemar Cantor to his cousin Dmitry Meyer in 1851: "Lelya, or namely Ludwig Gustav, the second, six-year-old boy, is possessed of an extraordinary glibness (frequently within logic), he has temperament of choleric person, and at times he is desperately stubborn".

In 1847 Georg Woldemar Cantor suffered from exacerbated tuberculosis[3], and he visited Italy for a course of treatment, along with his spouse. Their children were left for the trust of their neighbours, Anastasia and Osip (Ioseph) Grimm – the aunt of Georg Woldemar Cantor and her husband, a violinist.

According to the letter written by Georg Woldemar to his cousin, "After having a long travel within Southern Europe being accompanied by my spouse, in order to restrain the tuberculosis acquired through my labour (not inherited), after a long staying in Pisa and "belle Venice", in the autumn of 47$^{th}$ I arrived promptly in St. Petersburg, being more exhausted than before undertaking my voyage, for the only purpose of getting the blessing and leave this world in the presence of the two children[4]". However, his health, and then, the capacity for work, and well-being returned to him.

Sophia was born in the Cantors' family in 1848. By the end of 1849, the fourth child was born - Konstantín Karl[5]. From the cited letter by Georg Woldemar Cantor to his cousin: "The girl which is nearly 4 years old, named Sophia, is an affectionate, charming and at the same time intellectual creature. Konstantín Karl is the last our boy, who is 2 years old; his temper is not defined yet".

Afterwards, Sophia proved to be an excellent painter, and it was she who instigated the interest of grown-up Georg Cantor in studying the Shakespeare's creations. Konstantín Karl became a good pianist and later, cavalry captain.

In that same house, P.L. Chebyshev took up residence since the year of 1850, and later, his friend, O.I. Somov, university professors, mathematicians [2].

When going out from the house Chebyshev and Somov could see how Cantors' children were playing in the garden, and hear the music from their windows.

The wife of Georg Woldemar Cantor, Maria Cantor, nee Boehm, originated from a famous musical family. She was enchanting, affable and artistic, she played violin. Their house was open for guests. Their friends-musicians visited them, and the music was played. She taught her children to play violin and piano. Georg Cantor cherished the affection for violin music for all of his life, and during his studentship he even organised a string quartet [3, p. 278].

Later on, Georg Cantor was sorry that his father did not allow him to become a violinist [4, p. 278].

Maria's brother, Ludwig Boehm, a violinist, used to come to visit them, as well as Anastasia Grimm, Georg Woldemar's aunt, with her husband Osip, a violinist; Georg Woldemar's uncle, Gartwig Meyer, a violinist, who loved to smoke a Havana when visiting his nephew. His daughters Masha and Mila, cousins of Georg Woldemar, also used to visit them. Also visiting were Anna, the sister of Maria, and her aunt, Justina Moravek, both served as maids of honour[6] at court of the Grand-duchess Maria Nikolaevna. A violinist and composer Ludwig Maurer, a friend of the late father of Maria Cantor-Boehm also visited them. Charles Moberly, a broker, a friend of Georg Voldemar Cantor also came to visit them.

When Georg was eight, he and his brother were sent to study in St. Petri-Schule (Nevsky prospect, No. 22), the Principal German specialized school under St. Peter's Lutheran Church. (Picture 8. Petrischule). The special school was established in 1709, for numerous German speaking population of St. Petersburg. The building is preserved in the shape it had after rebuilding during 1830s. The school is located behind the building of the Church.

---

[3] At that time, no information was available about infectious nature of tuberculosis. The first work by N.I. Pirogov on the subject matter was published in 1852.

[4] Georg Woldemar Cantor was inaccurate in both the dates and indication of exact age of his children.

[5] Dauben wrote, "The couple`s six children, of which Georg was the eldest" [4, p. 274]. No information was found about other children in addition to the data on the four aforementioned children.

[6] The term "maid of honour" (*Kammer-jungle – German*) stands for a maiden who assists as dressing of royal persons. To form a staff of Grand-duchesses, young daughters of musicians and merchants able to participate in games and musical lessons were selected.

Teaching was of high level, many people, not only Germans, were seeking to send their children to the school. The school had 15 classes, with 20 to 70 people in each class; there were classes for boys and classes for girls. There were 67 pupils in the Cantor's class.

The subjects taught in the first class included: religion, Russian, German, French languages, arithmetic, calligraphy, geography, drawing; optionally – English and Latin[7]; in the second class – religion, Russian, German, French languages, arithmetic, calligraphy, drawing.

During the first term Georg Cantor was a very good student, he had excellent marks: religion – 4, Russian language – 5, French – 4, arithmetic – 5 [6]. However, during the second (spring) term of the year of 1853 everything turned to be very bad: religion – 3, 2, 2, Russian language – 3, 2, 2, arithmetic – 4, 3, 3, 2, other disciplines – 3, 2, 2, 2. The average score was "2.2" [7]. It was possibly due to his father's illness, or other family's problems, or maybe it was the first episode of depression (Picture 9. Classbook with Cantor`s marks).

In the Georg Cantor's class I. Laguzen was the teacher of penmanship (calligraphy), and Cantor's marks were 3.5 and 4; G. Lalans graduated from that same special school was the teacher of French, he gave Cantor marks from 3 to 4.5; C. Limmerich was the teacher of geography, and Cantor's marks were only "fives". Drawing was taught by V. Pape, and Cantor's marks were 4 and 5. German grammar was taught by R. Schultze, the marks were 4.5 and 5. D. Ulyanov was the teacher of Russian language. Cantor's marks in this course were from 3 to 4.5. Cantor did not learn Latin.

The teacher of mathematics in the first and second classes was F. Feichtner. He created a table of Russian coins, measures and weights for his students. Georg Cantor had "fives", rarely "fours" in mathematics.

In the list of students of the second class of the year of 1854: Georg Cantor, of nine years old; religion of the family – Lutheran and Catholic, from merchant's family, was born in Petersburg, lives with his family; his younger brother Ludwig Cantor, 8 years old, is in the first class [8].

In that years the family moved to Bolshaya Konyushennaya Street, into Zhadimirovskiy's four-storey building, nearer to the school [9]. Present address: Bolshaya Konyushennaya, No. 1.

In the third class Georg Cantor was taught mathematics by V. Zarnov. The course of mathematics included arithmetic, commerce and accountancy. Zarnov was known as compiler of tables for calculation of interest on bank notes [10]. Cantor's marks in mathematics became notably better [11].

By the end of the 1st term of 1856 Cantor had highest successes in geography, then in mathematics and German grammar, then in religion, French, calligraphy, and the least in Russian [11].

Unfortunately, tuberculosis of Georg Woldemar Cantor, Georg's father, made it impossible for him to stay in a severe climate of Petersburg. Besides that, financial crisis began in Russia after accession of Alexander II to the throne; conditions of trade were being impeded.

On May 4, 1856, Georg Voldemar Cantor asked for 1 year of leave to travel abroad: "In consideration of the request of stockbroker Georg Cantor, who needs to travel abroad for one year, the Stock Exchange Committee, in its own sight, not finding a reason for impeding dismissal of Cantor on leave, has the honour to inform the Department of Foreign Trade of the same" [12].

In 1856 the family of Cantors left Petersburg for Germany.

## Relatives on the mother's side

[7] Cantor did not learn Latin in St. Petri-Schule.

(Picture 10. Maria Cantor). The mother of Georg Cantor – Maria Boehm (Böhm) – was born in the family of Petersburg violinists, Franz Boehm (Franciscus Böhm) from Hungary and Maria Moravek (Maria Moraweck, Morawech, Morawek) from Vienna. Franz and Maria Moravek got acquainted and married in Petersburg. This happened as follows.

In a Hungarian city of Pest, in the family of Michaelis Boehm and Anna Dorfmeister (Michaelis Böhm & Anna nee Dorfmeister) [13], two sons were born, Franz in 1788, and Joseph in 1795. Both of them became violinists, however, not in their original homeland, but in other cities – Franz in Petersburg and Joseph in Vienna. Franz, Joseph and Maria Moravek were taught the violin by Pierre Rode [5, p. 278].

*Joseph Boehm.* (Picture 11 Joseph Boehm). Joseph became a soloist of Royal capella in Vienna, professor of Vienna conservatoire, he wrote pieces for violin. He is reputed to be the founder of Vienna violin school, among his students were such well-known violinists as J. Joachim, Ya. Dont, E. Rappoldi, E. Remenyi and his nephew from Petersburg, Ludwig Boehm, subsequently professor of Petersburg conservatoire. Joseph Boehm was close to F. Schubert and L. van Beethoven, took part in author's concert of Schubert; Beethoven entrusted him with performance of his last quartets [14].

*Franz Boehm*, the first violinist of Imperial theatres in Petersburg. (Picture 12. Franz Boehm).

Well-known violinist Pierre Rode served in Petersburg from 1803 till 1807, and during that time young Franz Boehm, Joseph Boehm and Maria Moravek were his students. Georg Cantor recalled in his letter to P. Tannery in 1896: "My grandfather and my grandmother, Franz and Maria Boehm, nee Moravek, from the school of Frenchman Rode in St. Petersburg, during 20s and 30s were imperial violinists-virtuosos and arose admiration of musical circles, and my great-uncle Joseph Boehm, also Rode’s student, established a well-known violin school in Vienna, among the graduates of which were Joachim, Ernst, Singer, Hellmesberger (father), L. Strauss[8], Rappoldi" [5, p. 278].

In 1809, Franciscus Boehm started to serve as a violinist in Imperial theatres of Petersburg. Director of the bureau of Imperial theatres, A.L. Naryshkin, a connoisseur of art and expert in people, discerned the twenty-year-old Hungarian’s talent and virile character. Franz was hired immediately as soloist [15].

*Maria Moravek.* Czech family of courtier maitre d’hôtel Leopold Moravek from Vienna [16] lived in Petersburg since 1787; among their eight children were their daughters Maria, born in 1795, and Sophia born in 1798. The family was well-to-do, the children were well-bred, and Maria was a talented violinist. She played so good that since the age of 17 she started to concertize. The war against Napoleon was over, and young Maria Moravek placed an ad in a newspaper: "Damsel Maria Moravek has the honour to inform that she will perform in vocal and instrumental concerto grosso on December 18, h. a., in the Philharmonic large hall, in which concerto she will play violin" [17], and she leased the large hall. Apparently, the concerto had a success, because next year she performed in two concertos, and in February, 1814, she performed together with Franciscus Boehm, and on July 8, 1814, she married him.

A record exists in the parish register of Saint Catherine’s Catholic Church in St. Petersburg: "Franciscus Boehm from Pest, Hungary, the son of Michaelis and Anna nee Dorfmeister makes a match with damsel Maria Moravek, the daughter of Leopold and Anna nee Mako Grosentes. Witnesses: Ferdinandus Gidello and Johannes Ronko" [13]. (Picture 13. Saint Catherine’s Catholic Church).

Neither her tender age, nor the birth of her children impeded Maria Boehm performing in yearly concertos during the fast. She performed together with her husband on May 3, 1822. In May of 1823 Maria died of breast disease. She was 28 years old. Franz Boehm was left with his four children – Adolf, Anna, Maria and Sophia. He continued with his serving as the first concerto performer in the Imperial theatres, he played in an orchestra and gave concertos [18].

[8] Ludwig Strauss (1835–1899) was born in Pressburg (Bratislava).

The second marriage of Franz Boehm was with Sophia, the Maria's sister, in 1824 (1798–1866) [19]. That marriage gave three children: Ludwig, Julia and Maximilian. (Picture 14. Ludwig Boehm).

Franz Boehm was the teacher of many well-known musicians, – among them were the composers M.I. Glinka, A.N. Verstovskiy, A.F. Lvov, and even the members of the tsar's family [20, p. 73]. Glinka wrote a solo for Boehm in 1836 in the opera "Ivan Susanin" [21, p. 272].

In the Petersburg publications there were many enthusiastic reviews on the concertos of Boehm. V.F. Odoevskiy loved his art very much, and called him the best of violinists, and he called the Boehm's fiddlestick "Silk bow of Cupid" [22, p. 107].

For over 30 years, Boehm was the first soloist of Petersburg. Franz Boehm died on February 16, at the age of 57, of the "suffering from his nerves' weakening". The notification about his death was placed in a newspaper "The first performer of concertos in the Imperial theatres, Mr. Franz Boehm, who belonged to the number of the most notable virtuosos of violin, had died lately, to the common regret of all the people who knew him, and the numerous admirers of his extraordinary gift" [23]. He was buried in Smolenskoe cemetery in the city of Petersburg.

In the Russian Museum, there is a painting by Gregory and Nikanor Chernetsovs, "Parade on the Tsaritsyn meadow", created in 1837. In the foreground there is a group of townspeople, who admire the parade. Those among the group are Pushkin, Zhukovskiy, Krylov, noblemen, actors, painters, musicians, 223 people altogether, - all the pick of Petersburg. To the right of the group, wherein Pushkin stands, in the group of musicians and actors, there is Franz Boehm (personage No. 185 according to the list). (Picture 15, 16 - "Parade on the Tsaritsyn meadow").

*Maria Cantor-Boehm.* The mother of Georg Cantor, Maria Boehm, the daughter of Franz and Maria, was a child of three years old when her mother died. Sophia Moravek, the sister of Maria Moravek, who became the wife of Franz one year later, committed herself to support the children. The family lived near to Theatre square, large families of actors and musicians lived around. Together, they arranged feasts and home concertos for children. Maria was affable, benevolent, she played violin. One day, a colleague of Franz Boehm visited their home, also a violinist of the Imperial orchestra, Lutheran, Dane, Gartvig Meyer, who lived not far from their residence. His nephew was with him, a young educated merchant from Copenhagen, Georg Woldemar Cantor. Maria and Georg seemed to like each other, and then they were married in 1842. She was Catholic, and he was Lutheran. What made him to come to Petersburg?

**Relatives on the father's side**

On January 6, 1896, Georg Cantor wrote to Paul Tannery: "My late farther Georg Woldemar Cantor, who died in Germany in 1863, was a child when he moved with his mother to St. Petersburg, and he was immediately christened as Lutheran. But he was born in Copenhagen (I don't know the year for sure, it was somewhere between 1810 and 1815), he was born from Jewish parents who belonged there to a Portugal Jewish community, and therefore, apparently, of Hispanic-Portugal origin" [5, p. 241]. (Picture 17. Georg Woldemar Cantor).

In XVII–XVIII centuries, merchants and specialists in the various areas used to willingly move to Petersburg. There was a need in military people, seamen, engineers, architects, medical doctors, musicians, painters, craftsmen. Many Danes came to stay here. They took up their residences in Vasilyevsky Island, in the lines that were nearer to the Stock Exchange [24].

*Abraham Meyer,* a merchant from Copenhagen.

Abraham Meyer, an elderly merchant from Copenhagen (died in 1801), great grandfather of the mathematician Georg Cantor, lived in St. Petersburg and was engaged in meeting the ships that delivered his merchandise to him.

In August, 1799, a Danish galiot named "Di Frau Maria" was at roadster near Kronstadt.

The ship delivered merchandise from numerous merchants – coffee, sugar, ginger, rum, wine, cardamom, cloves, textiles – altogether for a sum of 138 000 roubles.

On the left and on the right of the ship there was a plenty of free space for passage of vessels. A Russian naval vessel was departing from the Kronstadt jetty, and the wind moved the naval vessel towards the galeot. While the naval vessel was listing heavily, an anchor that hanged from the naval vessel struck against the galeot, and pierced a hole in the board of the galeot. Because of the fact that such huge anchors, having a weight of up to 4.5 tons, were hanged outside the boards of a naval vessel, protruding from her boards, such accidents were rather frequent. The galeot was holed and sunk rapidly.

The first go were the works on salvage and repairs of the galeot; the people did not suffer, however, all the salvaged merchandise have soaked. The soaked coffee started to give off smoke, the sugar was washed out, the wine escaped from the casks. There was nothing to be sold. The trades people (the vessel delivered merchandise from 20 merchants) have suffered losses, and Abraham Meyer has been in their number [25].

Somebody gave an advice to Abraham Meyer, to recourse to Paul I.

The elderly merchant had an audience with the tsar, and what did he ask? - To shape the destiny of his children! He had two sons – Osip and Gartwig, and two daughters. Then, in September, the young men were hired as court musicians, and on January 21, 1800, a decree was issued, "His Imperial Majesty has highly deigned to hire the elderly Copenhagen merchant Meyer Abraham, residing at present here in St. Petersburg, his children Osip and Gartwig Meyer, not employed anywhere, for the service in His Majesty's Theatres’ directorate as musicians, at an annual income of five hundred roubles for each of them".

In 1801, after the death of Paul I, Alexander I took the throne, and announced about maintaining all pensions and payments awarded by his father. He issued moneys monthly, against applications; all operations were made subject to rigorous accounting and reporting.

On June 26, 1801, the widow of Abraham Meyer, who had just died, visited the tsar. Obviously, she was in a great indigence and heavily in debt, because a decree of the department followed: "Deigning to concede the request of the Meyer’s widow, we most graciously order to give her two thousand roubles as a loan, with deducting the same amount during five years from the salary of her two sons, Kammer-musicians, Meyers. Alexander" [26].

*Gartwig Johann Meyer*, a court musician.

The two sons of the elderly Copenhagen merchant continued to serve as court musicians. Their elderly mother and one of their sisters still lived with them. Osip died in 1803. Gartwig Meyer, Lutheran, continued to serve as violinist in an orchestra. He married Charlotte Wulff, and they had children: Helena-Emilia, Natalia-Maria, Dmitry-Joseph, Alexander, Adolf. The Emperor Alexander I was the baptismal witness (formal godfather) of the Meyer’s son Alexander.

The Emperor often did such a favour to those served at court; it entailed the commitment to put a boy to service and provide a dowry for a girl. In fact, it was characteristic for the people who were close to the Emperor, and who were orthodox. A mysterious patronage, supported the family of Meyers, was of many years' standing, starting from the event of wreckage of "Di Frau Maria" galeot.

Gartwig Meyer served as a musician till the year of 1828. He died in 1867.

*Dmitry Ivanovich Meyer,* professor of law, the creator of Russian Civil Law. (Picture 18. Dmitry Meyer).

Dmitry Meyer (1819–1856) brought fame to the family of Johann Gartwig. He was called the farther of Russian Civil Law. He graduated from the Main Pedagogical Institute in Petersburg, the course of juridical sciences, with gold medal. He was sent to pass an in-depth training course in the University of Berlin, where he attended the lectures from 1842 till 1844. In 1845 he was appointed as an adjunct in the University of Kazan, and after defending his doctor’s dissertation "On the ancient Russian right of pledge" in 1848, as a professor. In 1853 he was elected as dean of the faculty of law. In 1855 he was transferred to the University of St. Petersburg. "It induces me to serve in Petersburg as in the centre of our intellectual life, and besides, the city with which I’m connected by blood relations" [27, p. 11]. Unfortunately, Dmitry

Ivanovich had only delivered a single lecture there. He died of tuberculosis on January 18, 1856. He was buried in Smolenskoe cemetery. (Families of Meyers and Boehms used to bury their relatives there).

Here are some memoirs about this remarkable man. "Students of the University of Kazan were given such a mass of knowledge, which the audience could not be able to obtain anywhere else during that epoch. In addition to extensive material arranged as a strict scientific system, the Meyer's lectures were imbued with such a human character, courage of feelings, which should have to affect the students in a fascinating manner. When the voice of protest against serfdom, functionaries' bribe-taking, the voice against differences in human rights by estates and creeds was heard from a rostrum during the 40s, one would have to come to a conclusion that the professor was possessed of a considerable civic courage. Audacious word of the teacher was not vanished without an effect on the students: An occurrence is known, where one of the Meyer's students refused to make a bargain of buying serfs, just because of the influence of his impression got out of the university [28, p. 34].

The writer Leo Tolstoy was Meyer's student: "When I was at the University of Kazan, I didn't actually make anything during the first year. The second year I began to study. At that time professor Meyer was there, and he took an interest in me, and gave me a work – comparison of the "Order" by Catherine again the "Esprit des lois" by Montesquieu. And I remember, I was captivated by the work, and I retired into the country (in the summer of 1846), I began to read Montesquieu; the reading opened infinite horizons before me; I began to read Rousseau and decided to quit the University, just because I wanted to study" [29, p. 148].

Meyer bequeathed his library to the University of St. Petersburg. The Meyer's works are still being republished now; these are being read and studied.

*Anastasia Grimm*, the daughter of Abraham Meyer, the aunt of Georg Woldemar Cantor.

Dauben refers to Georg Woldemar Cantor's aunt, who married Joseph (Osip) Grimm, a violinist, in Petersburg [4, p. 272]. That was her second marriage.

Initially, having arrived from Copenhagen, she adopted orthodoxy, under the name of Anastasia (Nastasya) Mikhailova, she got married, and then she became a widow. For the second time she married Osip Grimm, a violinist, in 1809. On November 3, 1809, their church wedding ceremony took place in Sergiyevsky cathedral. "Kapellmeister of His Excellency the Honourable Minister of War Count Alexey Andreevich Arakcheev, Osip Grimm being under Catholic law without changing the same, is married to baptized in holy christening from Jewish named Anastasia Michailova, betrothed under the order of St. Petersburg Ecclesiastical Consistory under the No. 3115 in Sergiyevsky of All Artillery Cathedral, on November 3, 1809, the bridegroom in the first, and the bride in the second marriage" [30].

Joseph and Anastasia lived in the Vasilyevskiy Island, near to the Cantors' house. Their marriage was barren, Anastasia's son from the first marriage departed from Petersburg.

Quite possible, that the Grimms' family took on to look after the young Georg Woldemar Cantor, the father of mathematician Georg Cantor. The arguments for this include the following facts: Georg Woldemar Cantor was a child when his mother took him to Petersburg, and he was christened right away into Lutheran faith. His farther, Jakob Cantor, possibly stayed in Copenhagen. No information is available about the role of his mother in the life of Georg Woldemar. Apparently, he became an orphan soon. But her sister, Anastasia Grimm, probably, took care of the boy. Probably, the boy was brought up in her family.

After his marriage, Georg Woldemar Cantor with his young wife took up his residence next to Grimms, in a house not far away from them. He gave his wife's and aunty Grimm's regards to his cousin Dmitry Meyer in the letter to him: "special regards from my wife and aunty Grimm". This childish name – aunty Grimm, ranking with his spouse, points to the fact that in the Cantor's family they had equal roles. In 1847, Georg Woldemar, when tuberculosis developed in him, went abroad together with his wife, and their two infants, a year-old baby Gustav and two-year-old Georg were left in Petersburg. So, who was to take care of them? – Most likely, the aunty Grimm.

The fact that Osip Grimm bequeathed to Georg Woldemar Cantor, "Villmanstrand and temporary Petersburg merchant, Yegor Yakovlev, the son of Cantor" [31], by his last will signed on January 30, 1854, also witnessed to their close relations.

The next day after signing the last will, on January 31, 1854, Osip Grimm died of dropsy, at the age of 68.

A petition is available, forwarded by Georg Woldemar Cantor ("Yegor Yakovlev Cantor, Villmanstrand and temporary Petersburg merchant of the third guild") written in the Russian language on May 26, 1854, and addressed to the Bureau of His Imperial Majesty, about a request of grant of the remaining sum of the salary of Osip Grimm.

The amount of the remaining pension of Grimm made 57 roubles and 39 kopecks in silver. In addition, Cantor got chattels personal for 350 roubles, pecuniary property for 33 roubles, 6 shares of fire insurance society, – 690 roubles each, – for 4140 roubles in silver; shares of life income insurance society, – 81 roubles in silver each, – for 648 roubles in silver; altogether, 5171 roubles in silver [31].

Anastasia Grimm died between 1851 and 1854. Both, she herself, and Joseph Grimm were the guardians to Georg Woldemar, and he, in turn, took care of them at their old ages.

*The mother of Georg Woldemar Cantor, the daughter of Abraham Meyer* (the name is unknown).

Supposedly, Georg Woldemar, being a child arrived with his mother, and probably, with his father, in St. Petersburg, was placed in evangelical Lutheran asylum. Dauben wrote about that with a reference to Grattan-Guinness [32, p. 348] in the report by the Danish genealogical society [33]: "Though the history of the family is vague, undoubtedly, the father of Cantor, Jakob, lived in Copenhagen. […] The family was deprived of everything during the military assault [of Copenhagen] and moved to Petersburg, where the mother of Georg Woldemar had relatives. Nevertheless, following the family`s move to S.Petersburg, the upbringing and education of the child Georg Woldemar was entrusted to an evangelical Lutheran mission. What became of his parents at this time is unknown, though his mather`s family, by name of Meyer, was reportedly a respect and successful family in St. Petersburg. A sister was married to Josef Grimm (a Roman Catholic), who was a chamber musician of the royal court, and a nephew was establishes as a professor of law at Kazan. He was apparently instrumental in promoting the legal apparatus involved in freeing the Russian serfs in 1861, and as the family liked to recall, Tolstoy had once been his student. As for Georg Woldemar`s parents, nothing more is known exept for the fact (mentioned in a report of the Danish Genealogical Institute) that Jacob Cantor still alive in 1841, when he sent congratulations to his son on the occasion of his engagement" [4, pp. 272–274]. All the hypotheses and possible speculations as regards the fate of this branch of the family are presented in detail in the book [34].

*Georg Woldemar Cantor*, merchant and stockbroker.

It is known that his mother was from Meyers's family, and that he was the cousin of Dmitry Meyer, the son of Gartwig [35]. Because of the fact that in 1837 he was mentioned in the book of addresses as the stockbroker Yegor Cantor [36], whereas a man could only become a stockbroker upon getting 30-year-old age; then one can assume that the year of his birth was 1807. However, later findings changed that assumption; according to the record of his marriage, he was born in 1814.

Prior to 1833, Cantor was not included in the list of merchants. The first information about him appeared in the "Commercial gazette" on November 22, 1834: "Cantor & Co. – delivery of hemp in amount of 1855 poods, 175 casks of hempseed oil; Yegor Cantor, delivery of potash, 31 casks". Hereinafter, in 1835, "Danish national, Yegor Cantor, enjoying the rights of local merchant of the second guild, delivered 6 casks of potash = 153 poods and 25 pounds". During the year of 1835, G. Cantor shipped merchandise for the sum of 223 396 roubles 50 kopecks, with delivery for 5 509 roubles 0 kopecks. During the year of 1836, Cantor & Company shipped merchandise for the sum of 12 200 roubles, with delivery for 149 784 roubles. During the year of 1837 Cantor & Company shipped merchandise for the sum of 114 700 roubles, with delivery for 60 394 roubles.

Since the year of 1837, there were no references in the "Commercial gazette" regarding the Cantor's commercial operations. Is seems that the most probable supposition was that he associated his capital with some other Danish company, possibly with Asmus Simonsen. (The sender's address on the Cantor's letter to D. Meyer read: "Messieurs Asmus Simonsen & Co., for Mr. Georg Cantor"). That trading house also had the name of Sarepta house, or Mustard house, and it was located in Konnogvardeysky lane, No. 4. Their commerce was typical for Petersburg port – they exported hemp, lard, potash, canvas, and imported mining machinery, spices, rum, silk, engravings, paper.

Asmus Simonsen was an agent of Evangelical society of the Hernhuters, a Lutheran sect of an especially strict direction, also named "Moravian Brethren", or "Bohemian Brethren".

In Petersburg, Georg Cantor-father lived in Vasilyevskiy Island, near to the Stock exchange. In 1837 his address was indicated [36] in the list of stockbrokers: Cantor Yegor, line 3, No. 21, Vas. part, block 145, and in the list of merchant offices of Cantor Georg, and the firm of Cantor & Co.: line 5, No. 7, Vas. part. In the school register of his son he was referred to as merchant [37]. Georg Woldemar Cantor was hereditary merchant (remember Abraham Meyer and his galeot "Di Frau Maria"). Since the age of 19 he positioned himself in Petersburg as merchant and stockbroker. We shall notice that he used a passport with an incorrect indication of his age (apparently, that was because of a 30-year-old age qualification requirement for commercial activity in Russia). At the age of 24 he had already stopped trading (probably, the reason was an imbalance in commercial operations described above).

Starting from 1835, Georg Woldemar Cantor signed up to become one of Villmanstrand merchants. A little town of Villmanstrand (now Lappeenranta, Finland) situated in 220 km from Petersburg was a part of the Grand Duchy of Finland within the Russian Empire. The Finnish merchants were granted considerable trade concessions, many of foreign merchants signed up to become Finnish merchants in order to get advantages of the offshore territory, though they were only "false burghers", who almost never visited the cities where they were registered. Thereby, Georg Woldemar Cantor became a citizen of the Russian Empire.

Since 1816, Georg Woldemar Cantor's maternal uncle, Gartwig Meyer, was a violinist in the orchestra of Imperial theatres, where Franciscus Boehm played the first fiddle. All actors and musicians lived near to the theatre. That time it was common to visit each other, with the own musical instruments, to play Schubert's quartets.

Perhaps it was one of such home concerts, where young people, Georg Woldemar Cantor, a well-educated merchant, and Maria Anna Boehm, a charming violinist, got acquainted.

The church wedding ceremony took place on April 21, 1842, probably, in Saint Catherine's Lutheran Church in the Bolshoy avenue of Vasilyevsky Island, No. 1. (Leonard Euler also visited this church, and probably, Georg Cantor was christened there).

Georg Woldemar Cantor was Lutheran, and Maria Boehm was Catholic. Therefore, on April 22, 1842, in Saint Catherine's Catholic Church (where Maria's parents were married in 1814) another wedding ceremony according to Catholic rite took place: "Georg Woldemar Cantor, merchant from Villmanstrand, Lutheran, 28 years old, and mademoiselle Maria Boehm, daughter of Franciscus Boehm, 22 years old" [38]. The witness on the side of the bridegroom was Charles Moberly, a broker; on the side of the bride– Ludwig Maurer, a violinist. A signature is also available on the document, which could be read, to some degree of confidence, as 'Jakob Cantor' (Georg Voldemar's father). Hence, we can establish the year of birth of Georg Woldemar Cantor – 1814. (The same year is written on his gravestone).

To judge from the composition of guests, the family of Maria Boehm was not happy with that marriage. Nobody from the bride's family attended the wedding. The only bridal party was a good soul, Ludwig Maurer, violinist and conductor, a friend of Franz Boehm. Franz Boehm, probably, disapproved the choice of his daughter, because the bridegroom was without means and from a different social environment (a nice person, but bad musician).

The marriage inspired Georg Woldemar with new power and self-reliance. Whereas in 1839 he was in the list of brokers unable to pay excise duty [39], in 1848 he was again elected broker [40], and his successfulness grew gradually; in 1854 his name was in the fifth place in the lists of brokers by the amount of excise duty paid [41]. Georg Woldemar even joined insurance society for widows; in order to provide support for his wife in case of his death[9]. In 1848 he was again elected as broker, and his successfulness grew gradually; in 1854 his name was in the fifth place in the lists of brokers [42].

Due to Grattan-Guinness, Maria Boehm adopted religion of her husband [32]. But in 1854 she was referred to as Catholic in a school register of her sons, Georg and Ludwig [43].

We shall cite some fragments of Georg Woldemar Cantor's letter [1] to his cousin, Dmitry Meyer, in October, 1851. Though both of them were of Danish ancestry, the letter was written in Gothic German cursive [10]. In the letter, Dmitry was referred to as Ossa, a pet name of Joseph, Dmitry's second name. Georg Woldemar Cantor was also in a decline, like Meyer. He recommended Dmitry to use a respirator, and mentioned his article about respirator, published in "Medicinisce Zeitung Russlands. St. Petersburg" [44]:

"My dear old Ossa,

For hundreds of times I conversed with you in my thoughts on the various subjects, and each time I had a desire and need to right to you, but the reality of urgent business matters, which, in severe hurry-scurry of my life, demand giving his best from a man who exists due to commerce, and besides, the most ingenious competition that exhausts not only physically, but morally as well, however incredibly it could sound, keep me in a state of apathy without letting me get free since the time of your departure, which means a poor incentive to be in correspondence with my friends or my family.

Looking at me from the side, it may seem to everybody that the highs and lows in my business do not actually affect ne anymore, and apparently, the material mode of life outshines everything, but believe me, whereas the pragmatism and forced striving for getting profits must capture fully, as it may seem to be, a man having four growing children, the man, however, still has his living soul, inborn aspiration for beauty and good; in secular cares, he simply very rarely feels an urge to confidence and finds a reason for the same!! […]

A few weeks back your father informed me of the content of your latest letter, wherein you communicated information about the testamentation, which indeed, completely frustrated the dear, kind old man, and plunged him into melancholy.

My dear friend, tell me, was it actually such a need to disturb your old parents by declarations and tidings of that kind??

Is it true that everything is so dangerous about your breast? Did your corporal constitution suffer due to your pertussis that persisted for a long time, or do you feel that you have got rid of this corporal ailment? Write me about that, it's my earnest request.

Masha Molayms has settled perfectly well, and Mila Kalinin[11] has got accustomed completely; she's got two sweetest children, her elder girl is truly nice, affectionate and clever one. Masha's Sonya is also a nice girl. Your mother has been in poor health lately, and tormented greatly by a violent stomach ache. Bothe sisters used to be night nurses for her; and now she has got better, and she has risen to her feet again, whereas an excellent physician, doctor Benbek has treated her; and now he has become my family physician, as well. Your old father also does not look young, however, he is possessed of an excellent appetite, and from time to time he uses to smoke a real Havana at our place, with great enjoyment, but I would like you to see that all with your own eyes.

---

[9] Due to that, she received a pension from Russia in amount of 190 roubles after her husband's death.

[10] To translate the letter, it was initially necessary to turn to a German company in Hamburg, engaged in deciphering illegible texts.

[11] Sisters of Dmitry Meyer, Natalia Maria and Helena-Emilia.

In conclusion, for my own purely selfish reasons, I would like to stay here for some more time, in order to avail myself of your views, experience, knowledge, methods and advice; and in order to bring up my two wonderful boys properly, from the very beginning […].

Well, it's all for today, - my own bookkeeping, which needs to be updated all the time, and a late hour, after midnight already, prompt me that I should finish writing, and probably you are bored by reading such a lengthy message from Georg to Dmitry Kazansky; special regards from my spouse, and aunt Grimm[12]; in my thoughts, I whole-heartedly shake your hand. Your old Georg". The sender's address was: Messieurs Asmus Simonsen & Co., for Mr. Georg Cantor in St. Petersburg.

In 1856 Cantor's family left Petersburg for Germany. That was due to changes in economic situation in Russia, and due to illness of Georg Woldemar Cantor, tuberculosis. Probably, Georg Woldemar was stunned by a premature death of tuberculosis of his beloved cousin, Dmitry Meyer. In the newspaper "St. Peterburgische Zeitung" dated 15(17) of May, 1856, in the list of those leaving for abroad: "Georg Cantor, Villmanstrand merchant, with his wife Maria Cantor, and their under-age children: Georg, Ludwig, Konstantin and Sophia, with Laura Zundshroem, Swedish citizen, residing at the address of Bolshaya Konyushennaya, in Zhadimirovsky's house, No. 1".

Georg Woldemar did not plan to depart forever; - he took a leave for a year. However, their removal proved to be final.

Initially, they took up residence in Frankfort on the Main. By that time the family was reach already. In 1871, after reunification of Germany and issuance of Reichsmarks (Deutschemarks), their fortune added up to more than half a million of Reichsmarks.

By that time the family was reach already. In 1871, after reunification of Germany and issuance of Reichsmarks (Deutschemarks), their fortune added up to more than half a million[13] [45]. Georg's farther wanted him to be an engineer, however, the choice of the son was mathematics, and he entered the University of Zürich to study the discipline. Georg Woldemar Cantor died in 1863, in Heidelberg. Maria Cantor died in 1896, in Berlin.

Cantor recalled his childhood in Petersburg with warmth. In 1894 he wrote in one of his letters: "My first wonderful 11 years spent in that beautiful city on the Neva River, unfortunately, would never recur" [46, p. 16]. At the end of his life, Georg Cantor, being a son of Russian citizen, thought about a possibility of his foreign service in Russia.

In this way, a family appeared in Petersburg, the members of which came from Hungary, Czechia, Denmark. Members of the family possessed of the most brilliant talents, like the violinists Boehms – the brothers Franciscus, Joseph, and the son of Franciscus, Ludwig; young violinist Maria Moravek, court musician, violinist Gartwig Meyer, Dmitry Meyer, the creator of Russian Civil Law, merchants and tradesmen Abraham Meyer and his grandson, Georg Woldemar Cantor. Life of the family was closely connected with the cultural life of Petersburg, and Georg Cantor inherited this culture[14].

**Dissemination of the ideas of Cantor in Russia**

Cantor wrote his principal works of the theory of sets during 1880s. Mathematicians from Russia, who have visited the universities of Berlin, Gottingen, read the 'Crelle Journal' (then all the universities were subscribers to the same), got acquainted with the ideas of the theory of sets. The first publication on the Cantor's theory belonged to Pavel Florensky. Pavel Florensky took a great interest in the Cantor's theory of sets; yet in 1900, being a first year student, Florensky wrote a foreword to the dissertation "The idea of discontinuity as an element of world outlook" [47], wherein he affirmed the unity of views of mathematical ideas of Bugaev and Cantor as regards world outlook. Since the year of 1902, Florensky delivered reports at the sessions of students' mathematical circle of the Moscow mathematical society, where he introduced an idea of association of mathematics with philosophy: "I have no doubt that our teachers, at least some of them, - and among them I cannot fail to mention N.V. Bugaev, - were possessed, to a considerable extent, of such entire world outlook, in the centre of which mathematics was placed" [48]. In 1904 he wrote a work named "On the symbols of infinity" [49]. In that work, for the first time ever in the Russian literature, he retold the doctrine of Cantor holistically, while making an attempt to apply the theory to the world outlook. His enthusiasm promoted the creation of an atmosphere, auspicious for the development of the theory of sets. Since 1907, the ideas of Cantor were included in the course of lectures delivered by I.I. Zhegalkin in the University of Moscow.

---

[12] Anastasia Mikhailova, married name Grimm, the aunt of Dmitry and Georg Woldemar.

[13] If convert these monies into Russian silver roubles, the fortune of Cantors made approximately the price of the merchandise sunk with the galiot named "Die Frau Maria". (On the principle that silver rouble contained 18 grams of silver, and Reichsmark contained 5 grams, the result makes amount of the fortune in roubles – 138,889 roubles – not a big sum of money for Petersburg merchants). To the point, we should take notice that posthumous debt of A. Pushkin made 140,000 roubles (100,000 – to private persons and 40,000 – to the Tsar).

[14] At greater length look in Russian [34].

# Illustrations

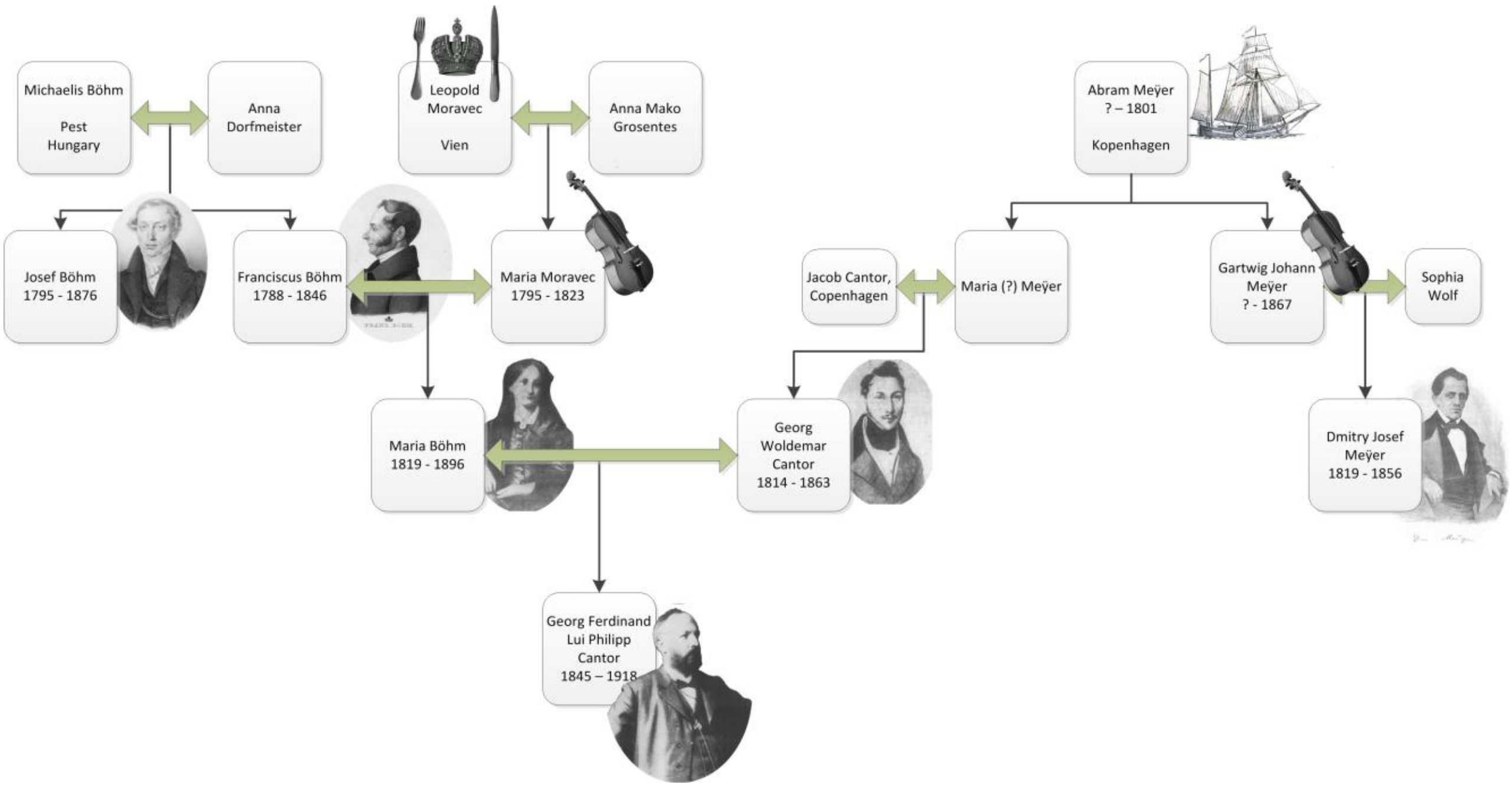

Picture 2. Fragment of G.W. Cantor’s letter.

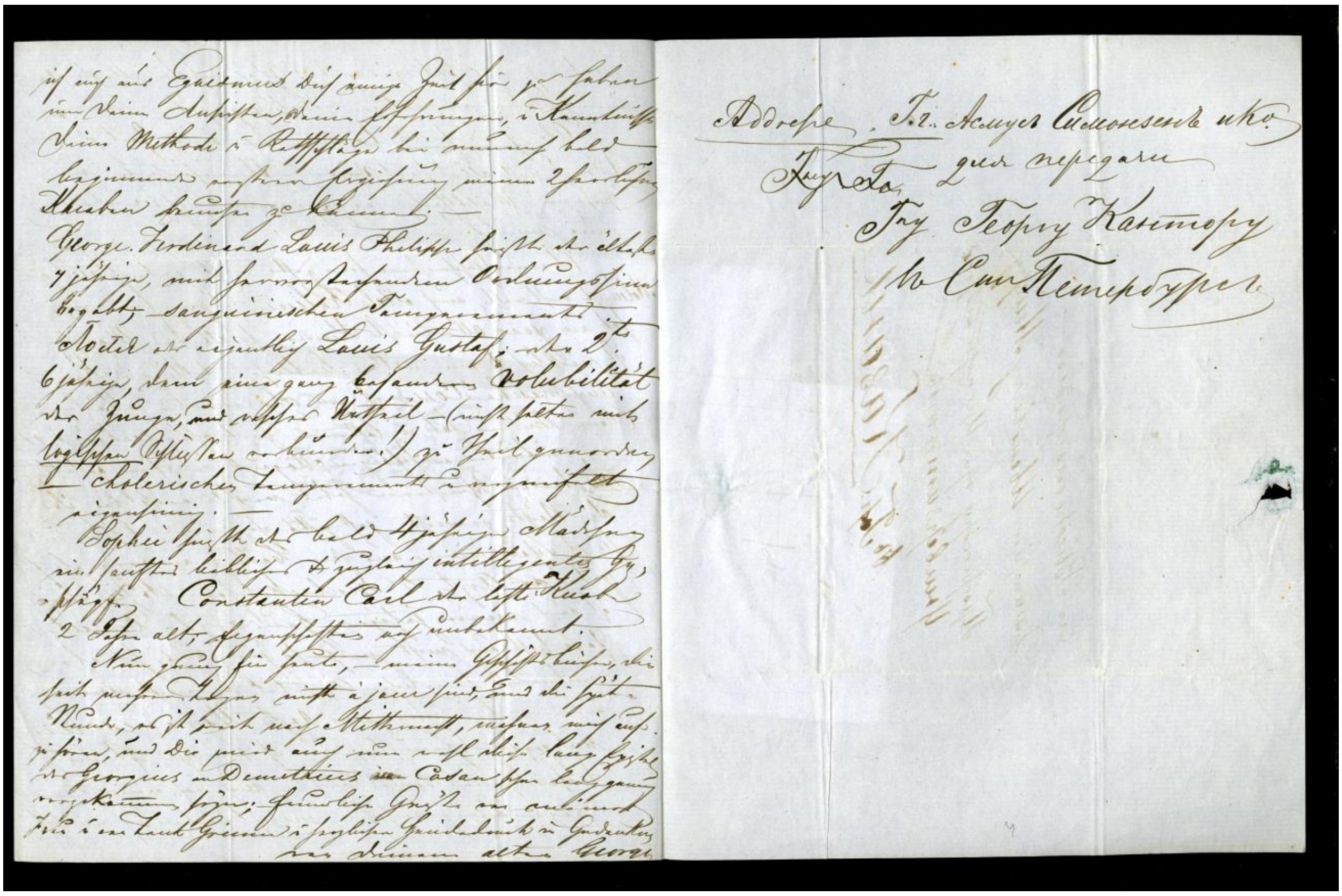

Адресъ. Гг. Асмусъ Симонсенъ и Ко.
Контора для передачи
Гну Георгу Кантору
въ Ст Петербургъ

3. House where Georg Cantor was born

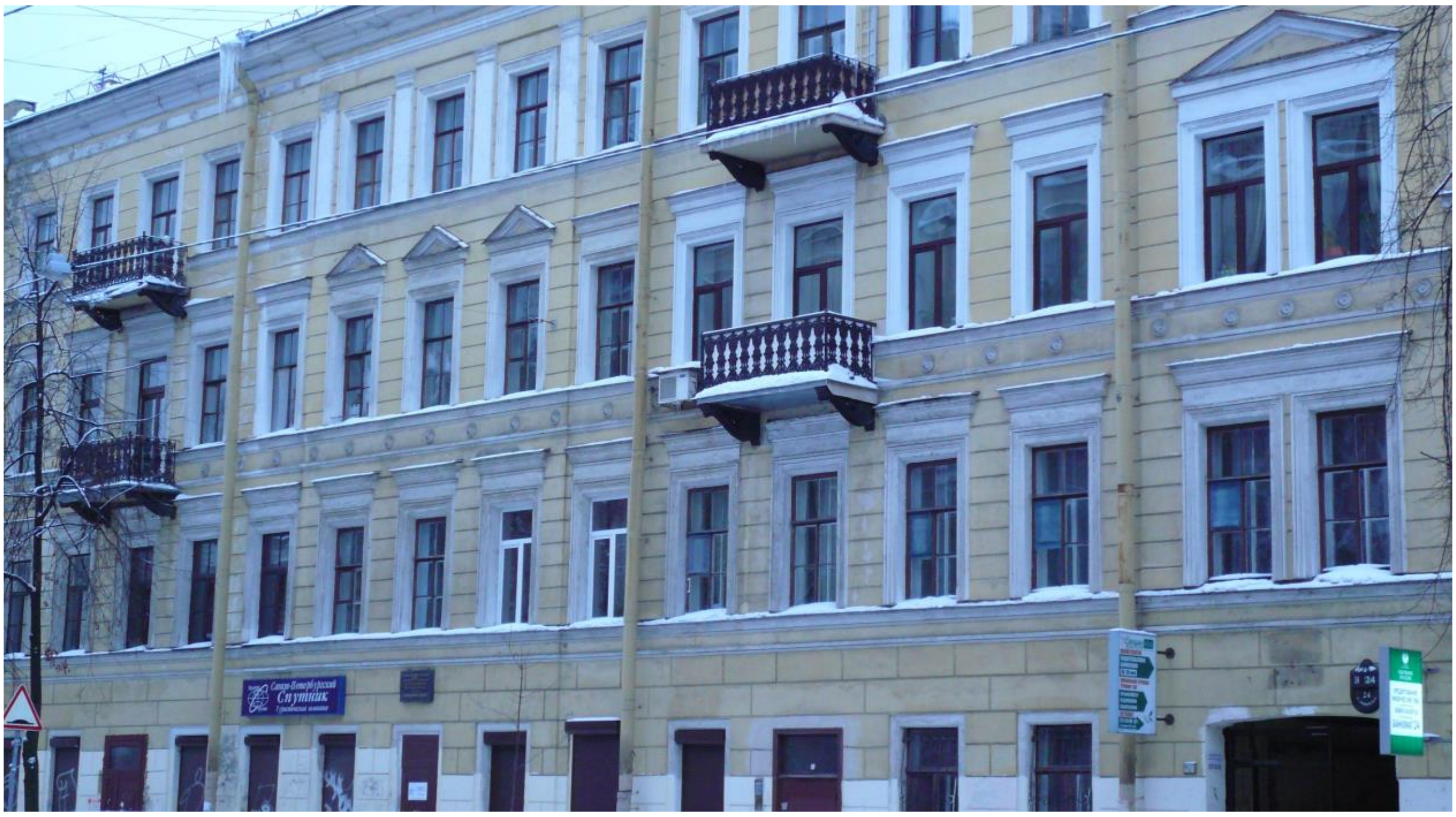

4. House where Georg Cantor was born

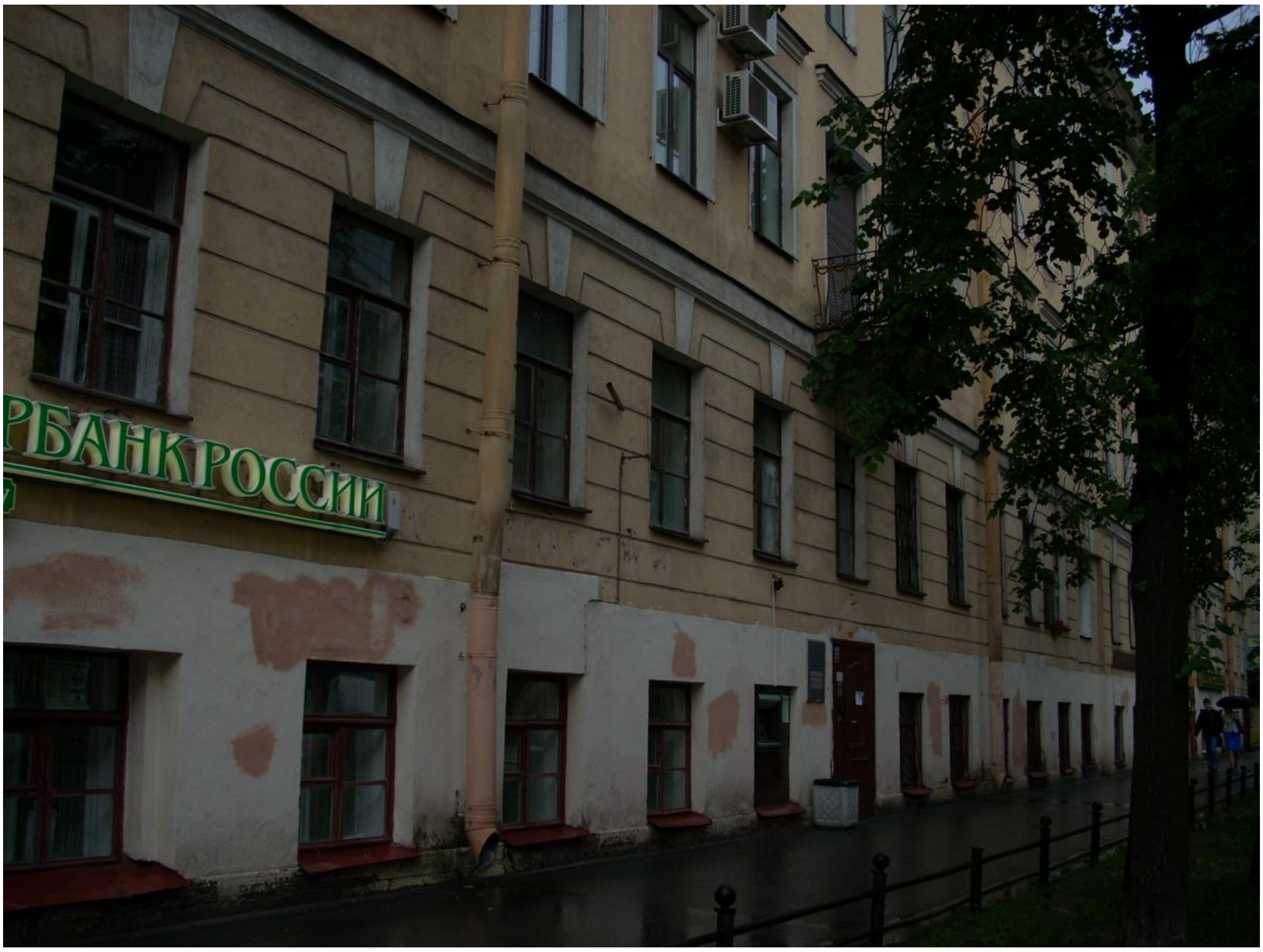

5. A memorial plague.

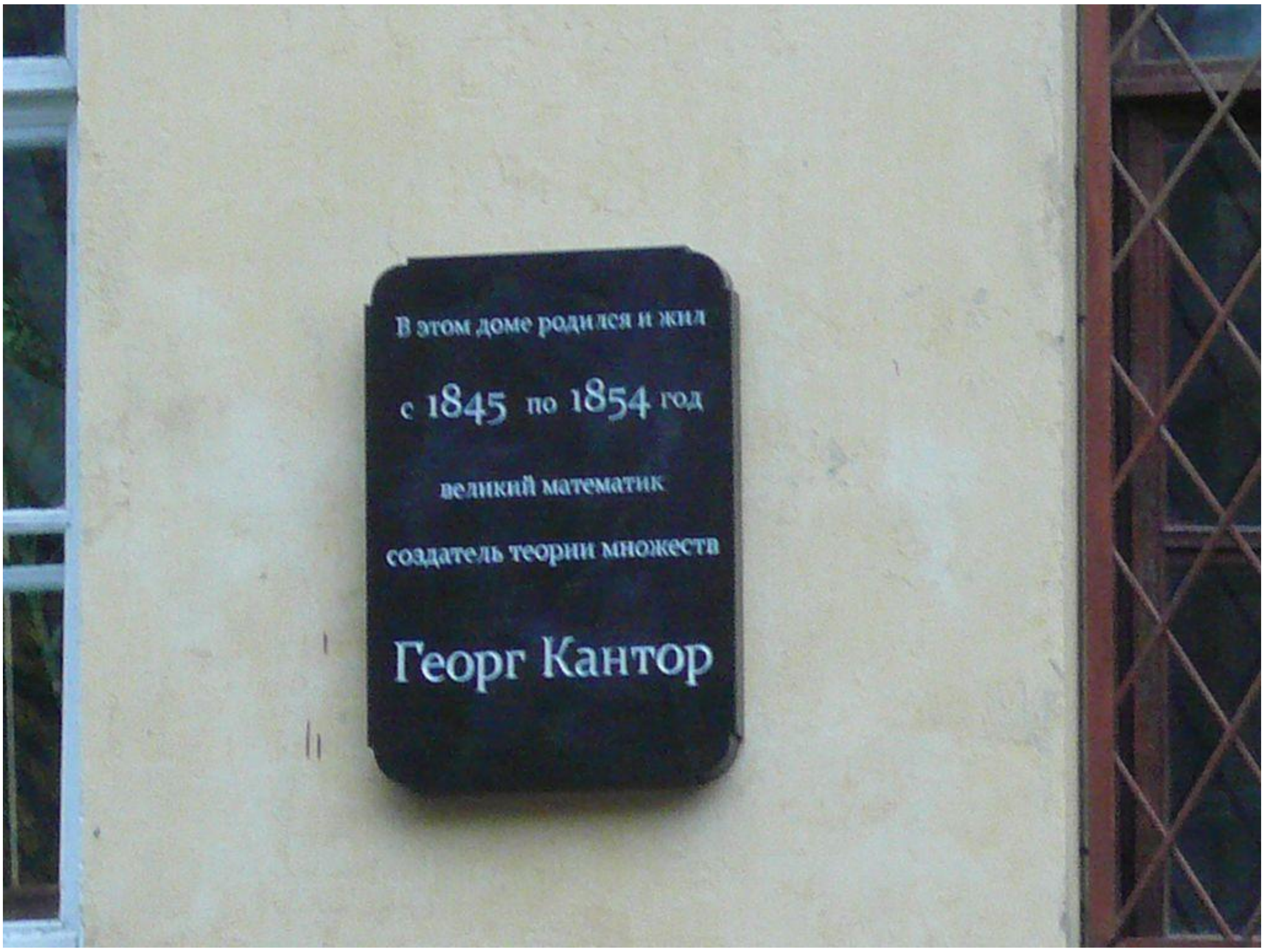

6. Stock Exchange

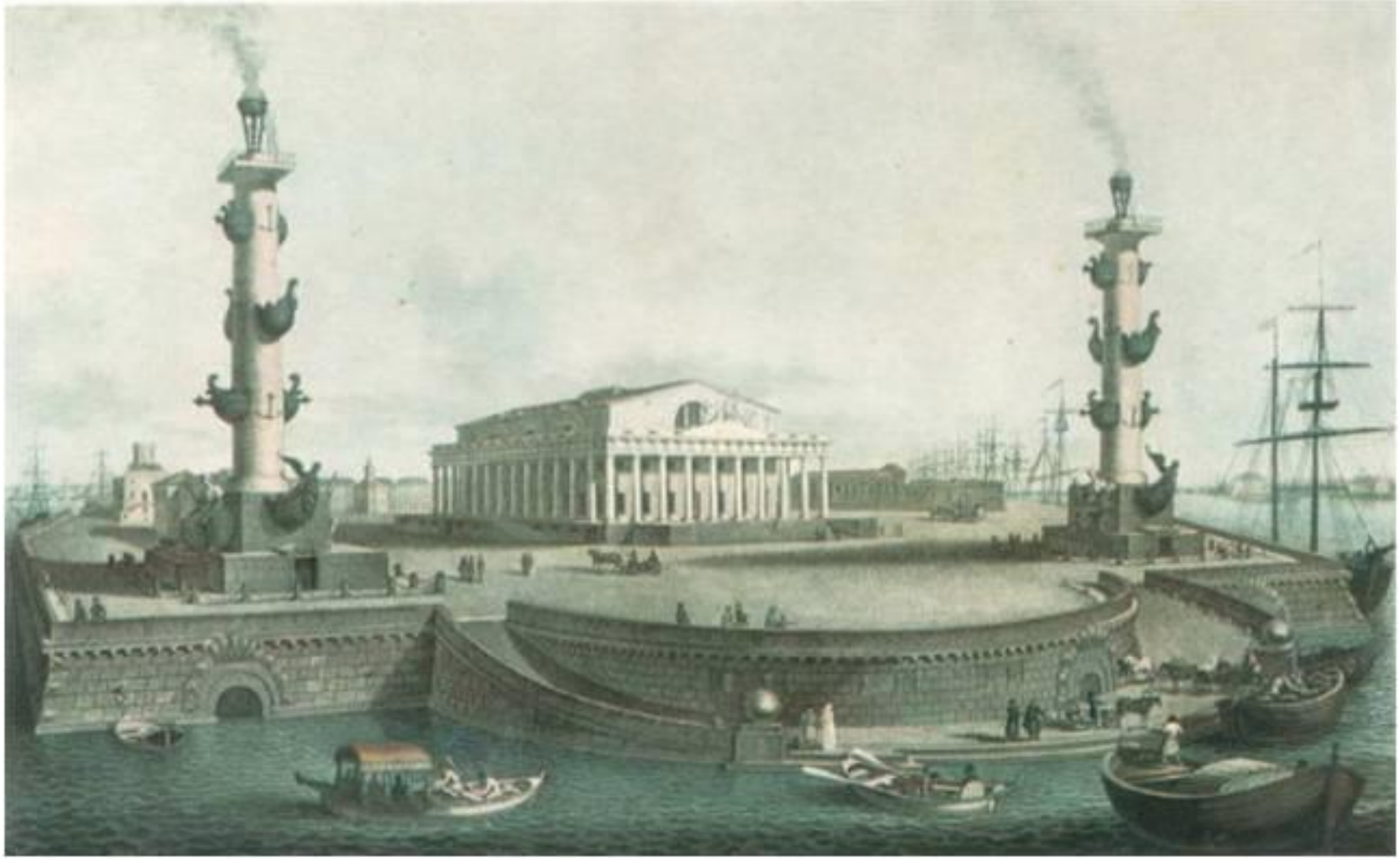

Picture 7. St. Catherine Lutheran Church.

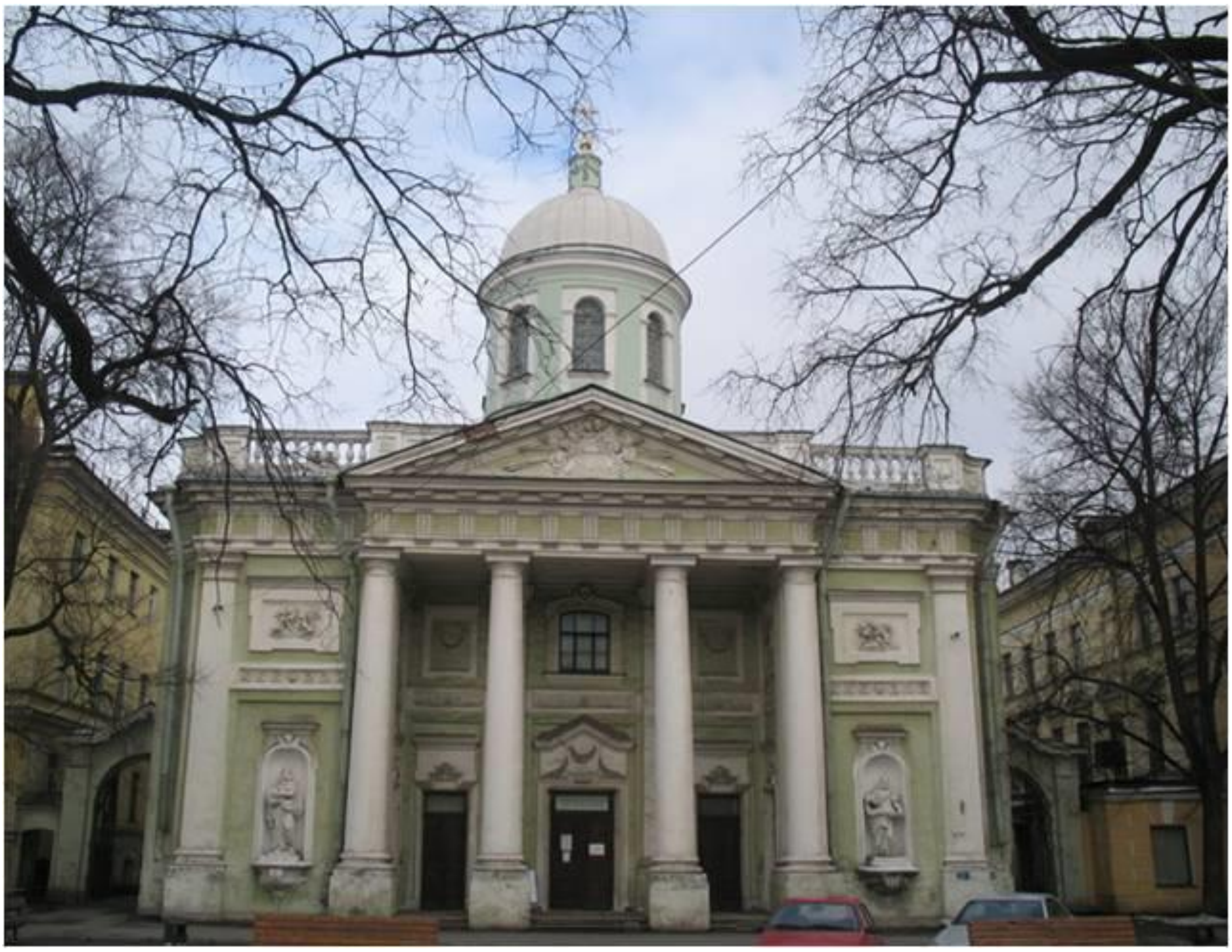

Picture 8. Petrischule.

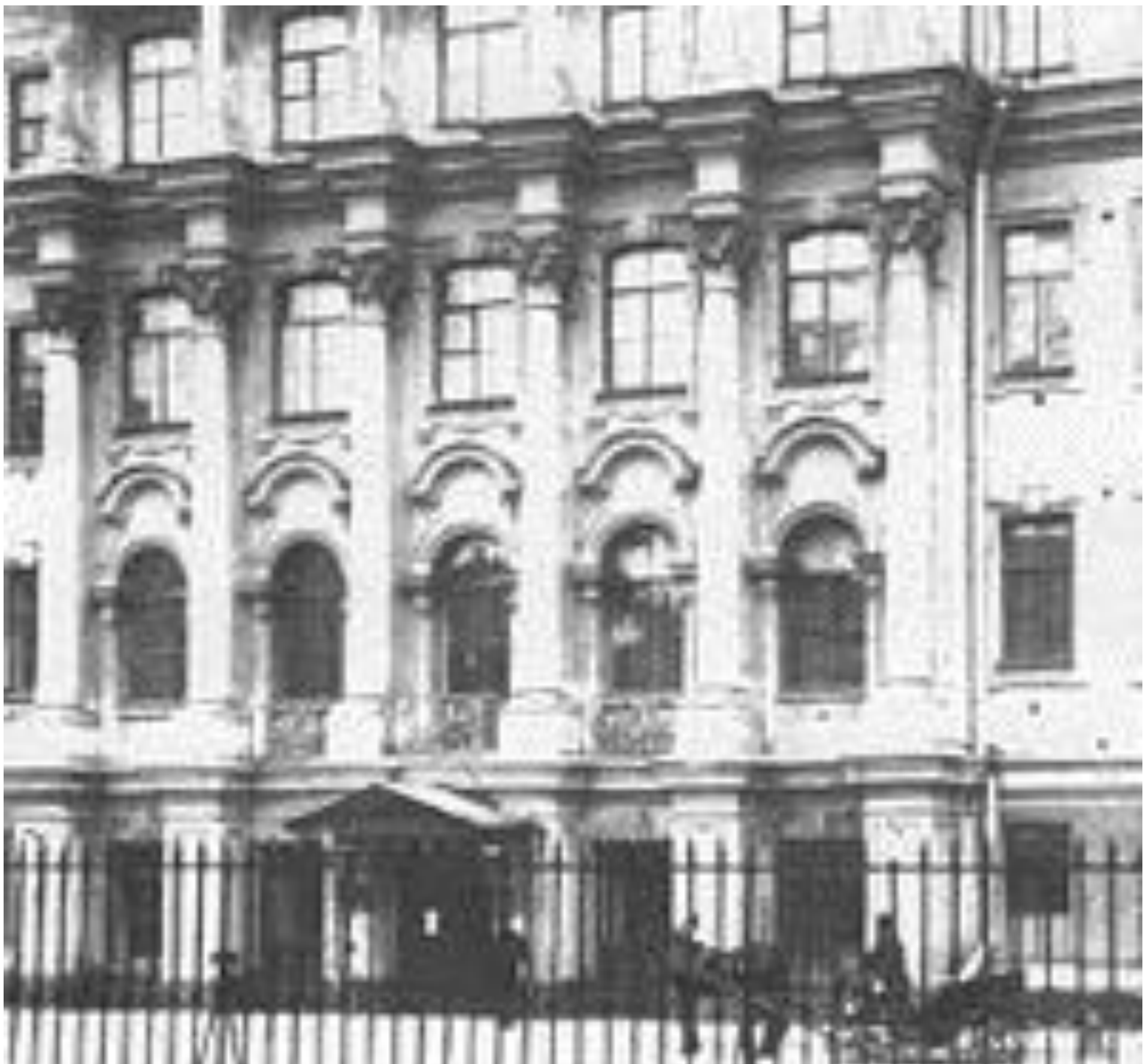

Picture 9. Classbook with Cantor`s marks

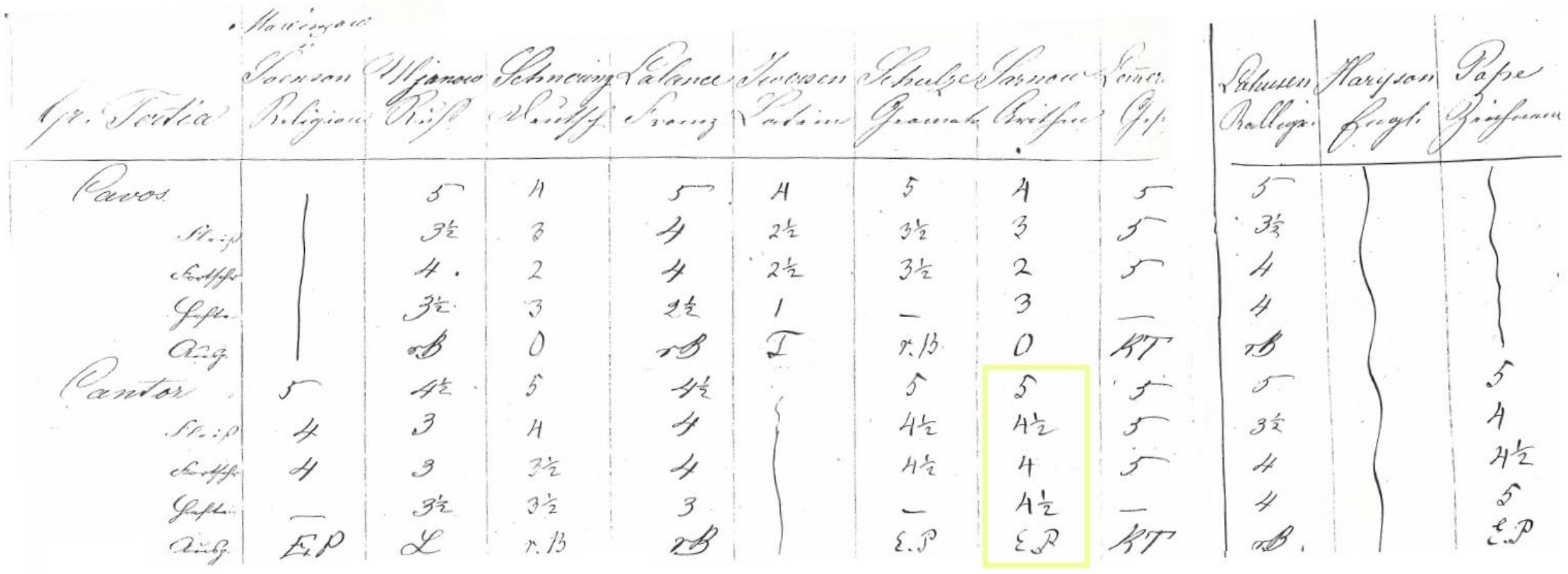

| Gr. Tertia | Jarnson Religion | [illegible] | Schnewing [illegible] | Lalance Franz | Iversen Latein | Schultze Gramm. | Jarnow [illegible] | [illegible] Ges. | Petersen Kalligr. | Harryson Engl. | Pape Zeichnen |
|---|---|---|---|---|---|---|---|---|---|---|---|
| Pavos | | 5 | 4 | 5 | 4 | 5 | 4 | 5 | 5 | | |
| Fleiß | | 3½ | 3 | 4 | 2½ | 3½ | 3 | 5 | 3½ | | |
| [illegible] | | 4. | 2 | 4 | 2½ | 3½ | 2 | 5 | 4 | | |
| [illegible] | | 3½ | 3 | 2½ | 1 | — | 3 | — | 4 | | |
| Aufg. | | rB | 0 | rB | I | r.B | 0 | KT | rB | | |
| Cantor | 5 | 4½ | 5 | 4½ | | 5 | 5 | 5 | 5 | | 5 |
| Fleiß | 4 | 3 | 4 | 4 | | 4½ | 4½ | 5 | 3½ | | 4 |
| [illegible] | 4 | 3 | 3½ | 4 | | 4½ | 4 | 5 | 4 | | 4½ |
| [illegible] | — | 3½ | 3½ | 3 | | — | 4½ | — | 4 | | 5 |
| Aufg. | E.P | L | r.B | rB | | E.P | E.P | KT | rB. | | E.P |

Picture 10. Maria Cantor

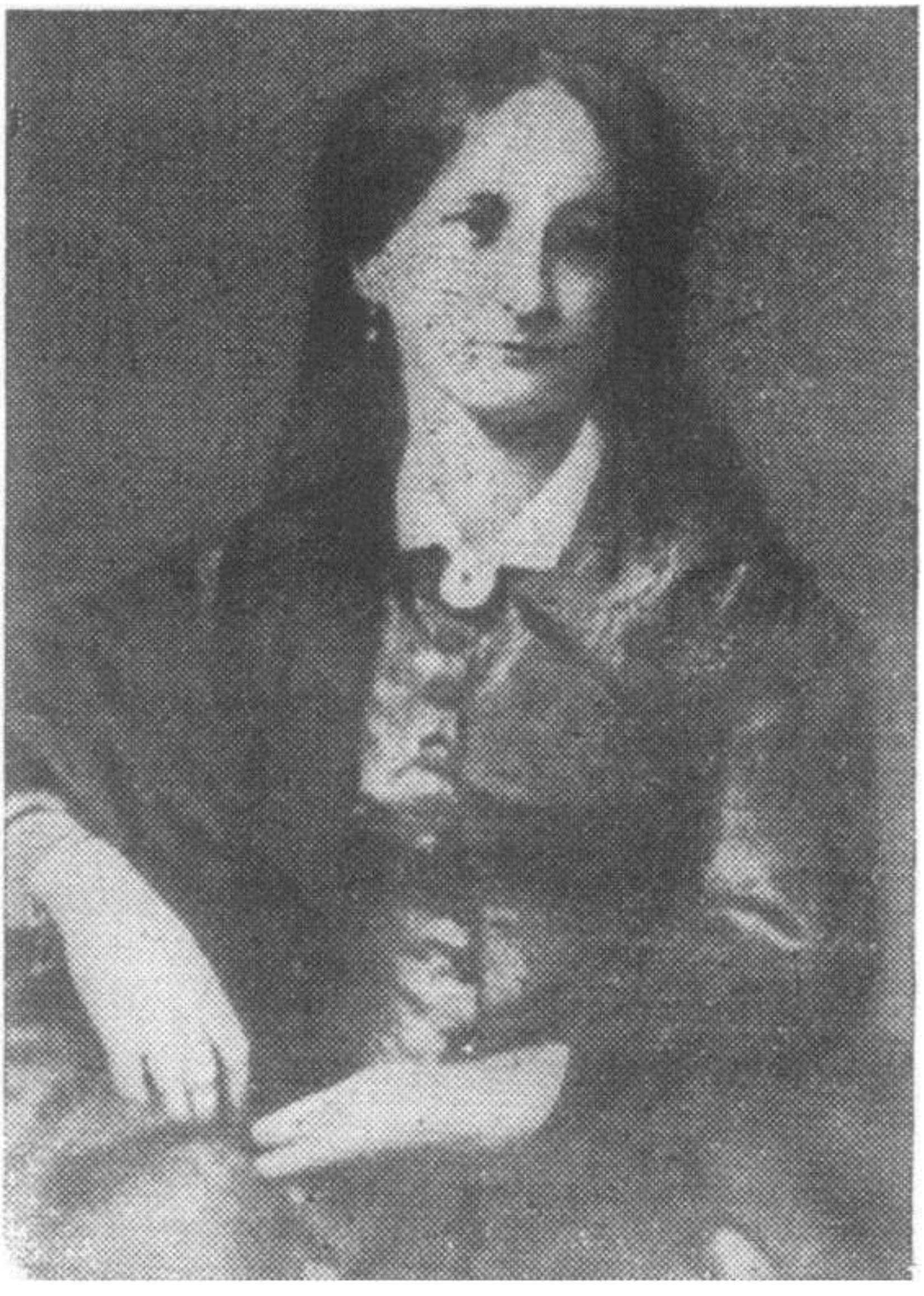

Picture 11 Joseph Boehm

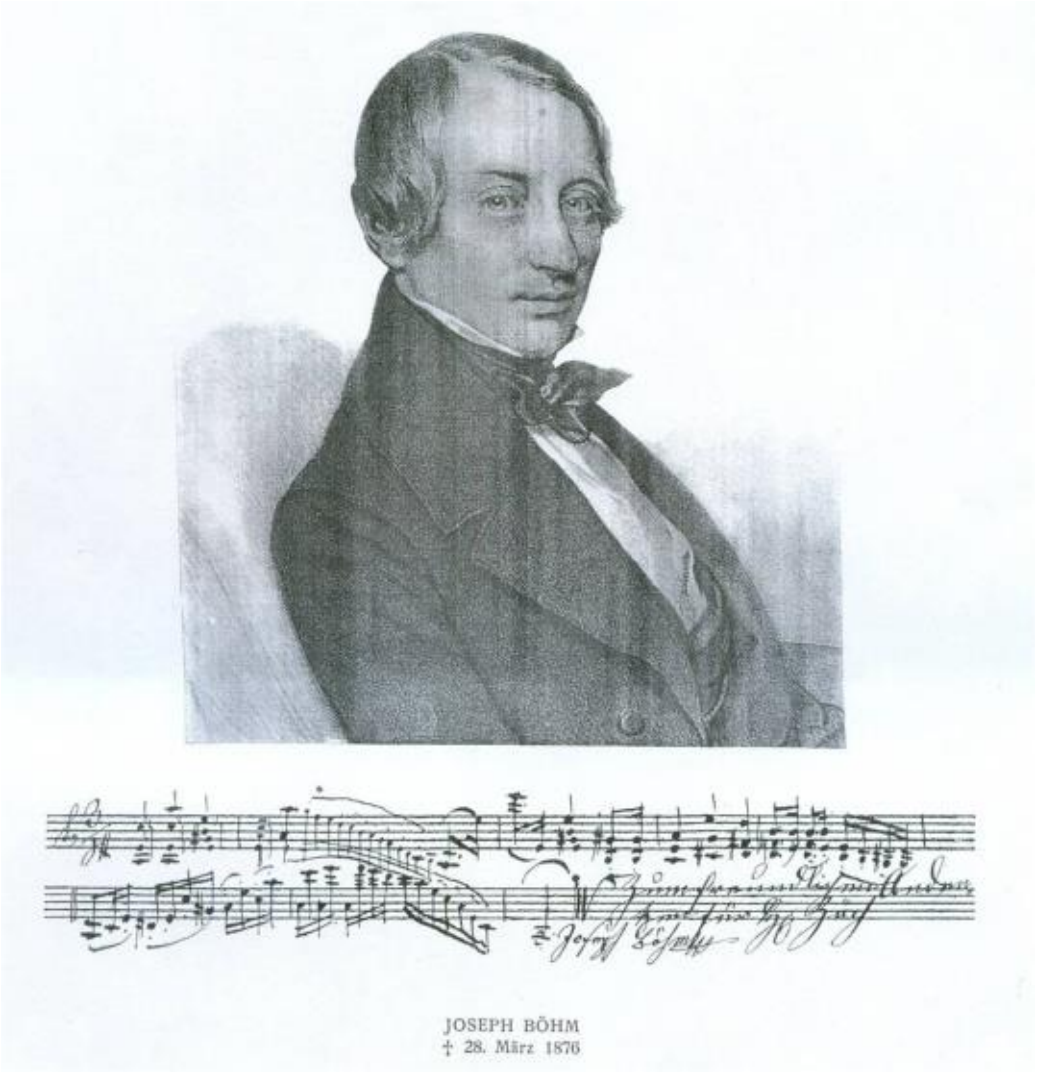

Picture 12. Franz Boehm

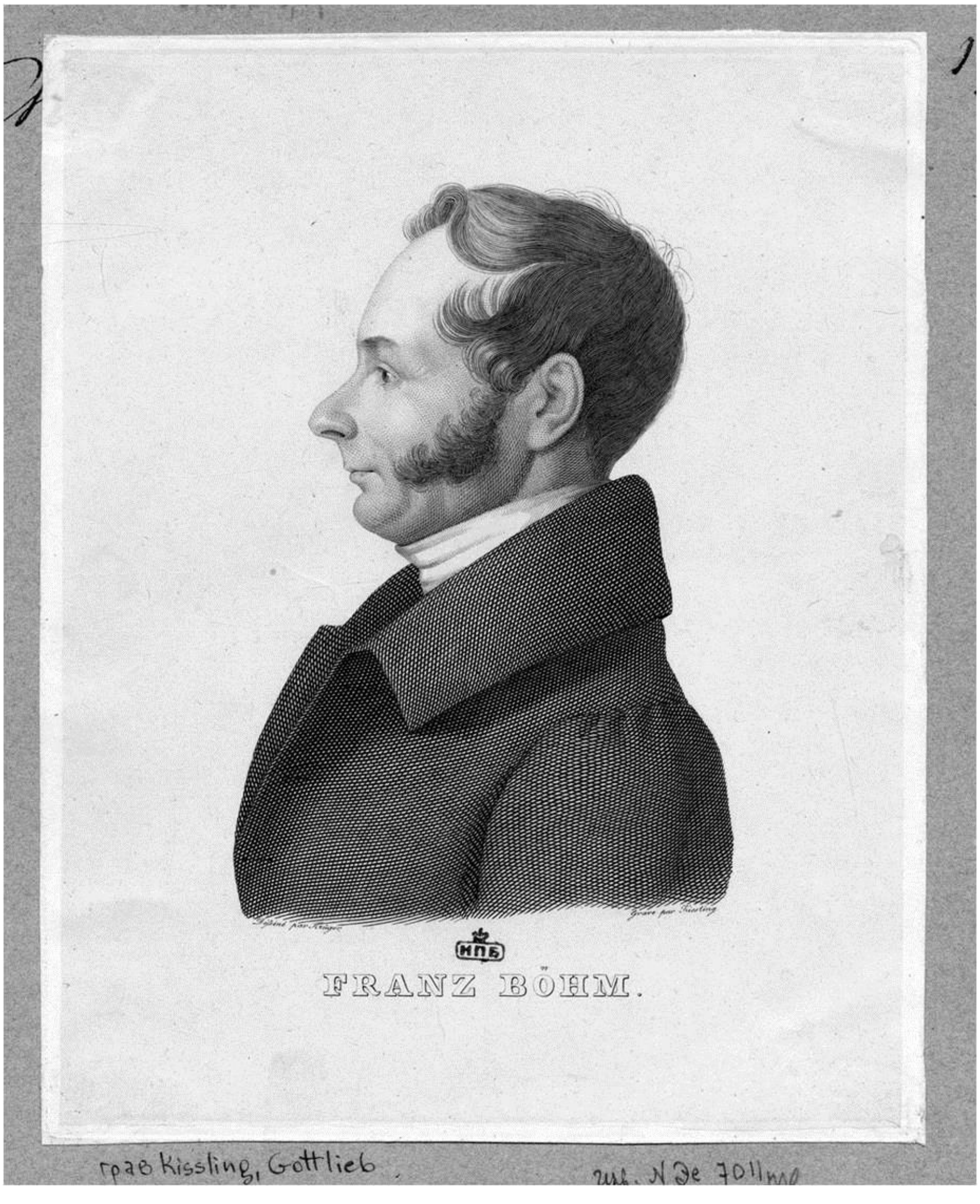

Picture 13. Saint Catherine’s Catholic Church

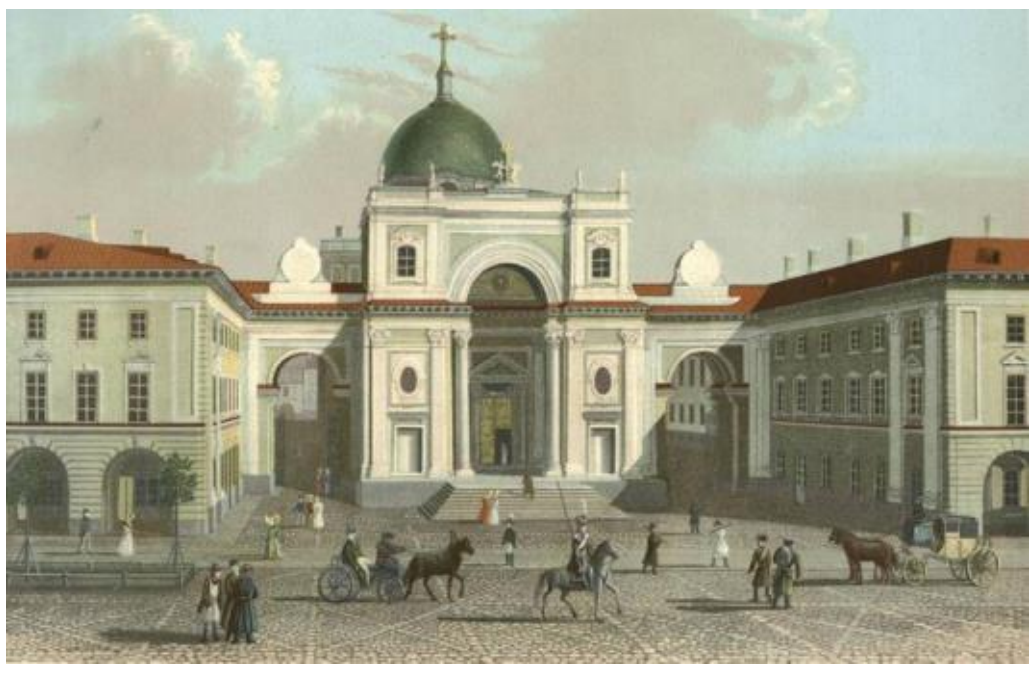

Picture 14. Ludwig Boehm

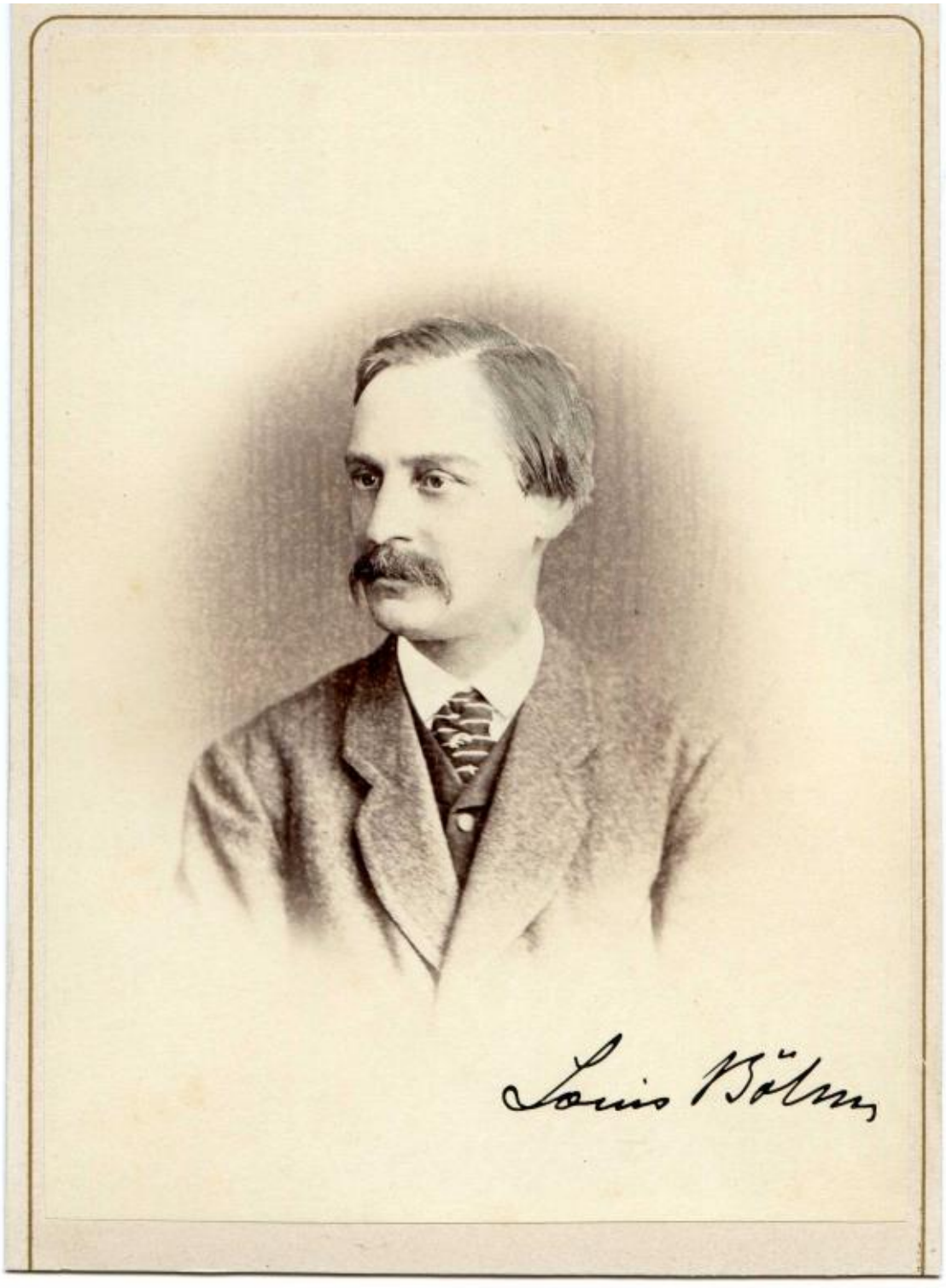

Picture 15, 16 - "Parade on the Tsaritsyn meadow”.

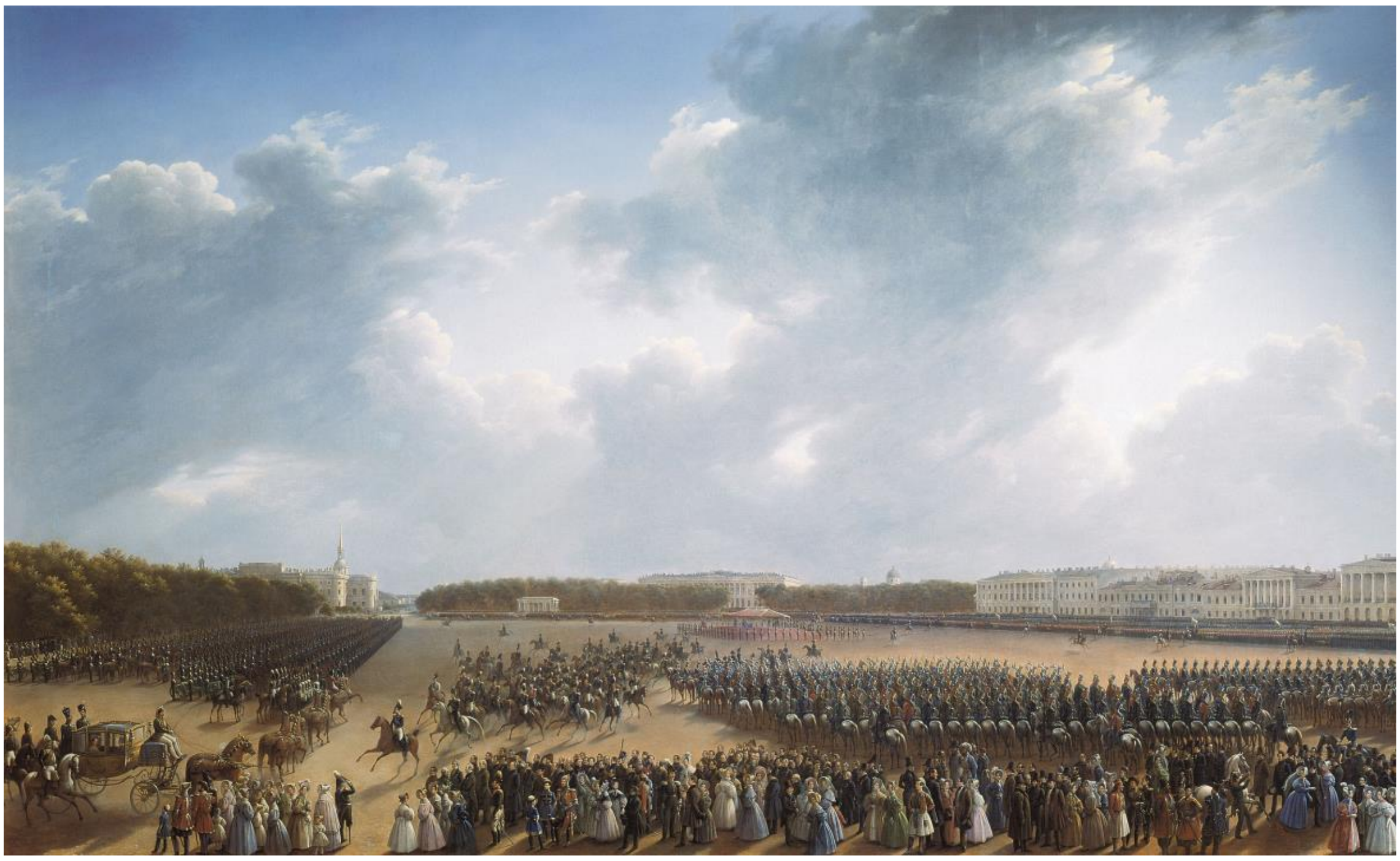

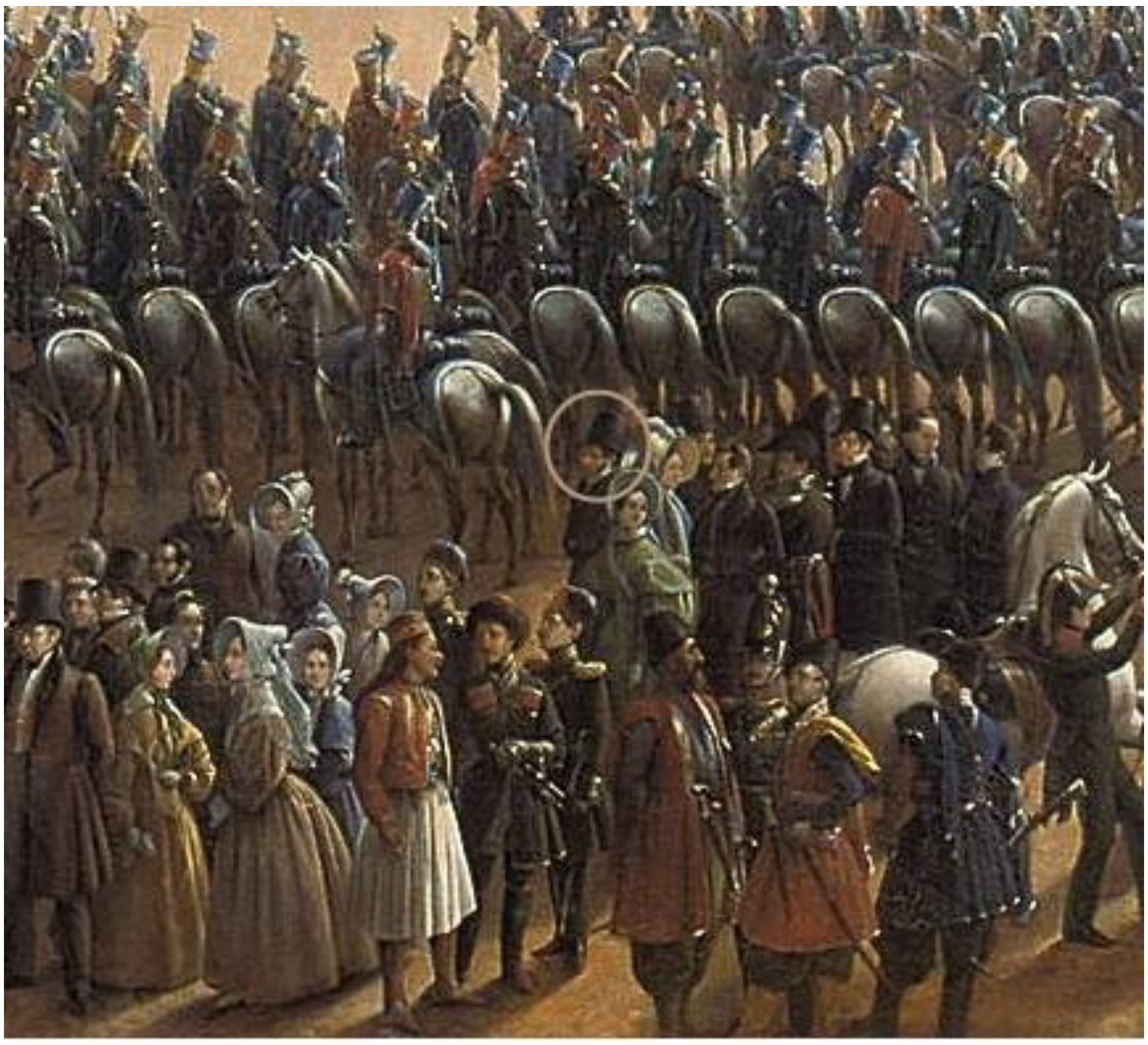

Picture 17. Georg Woldemar Cantor.

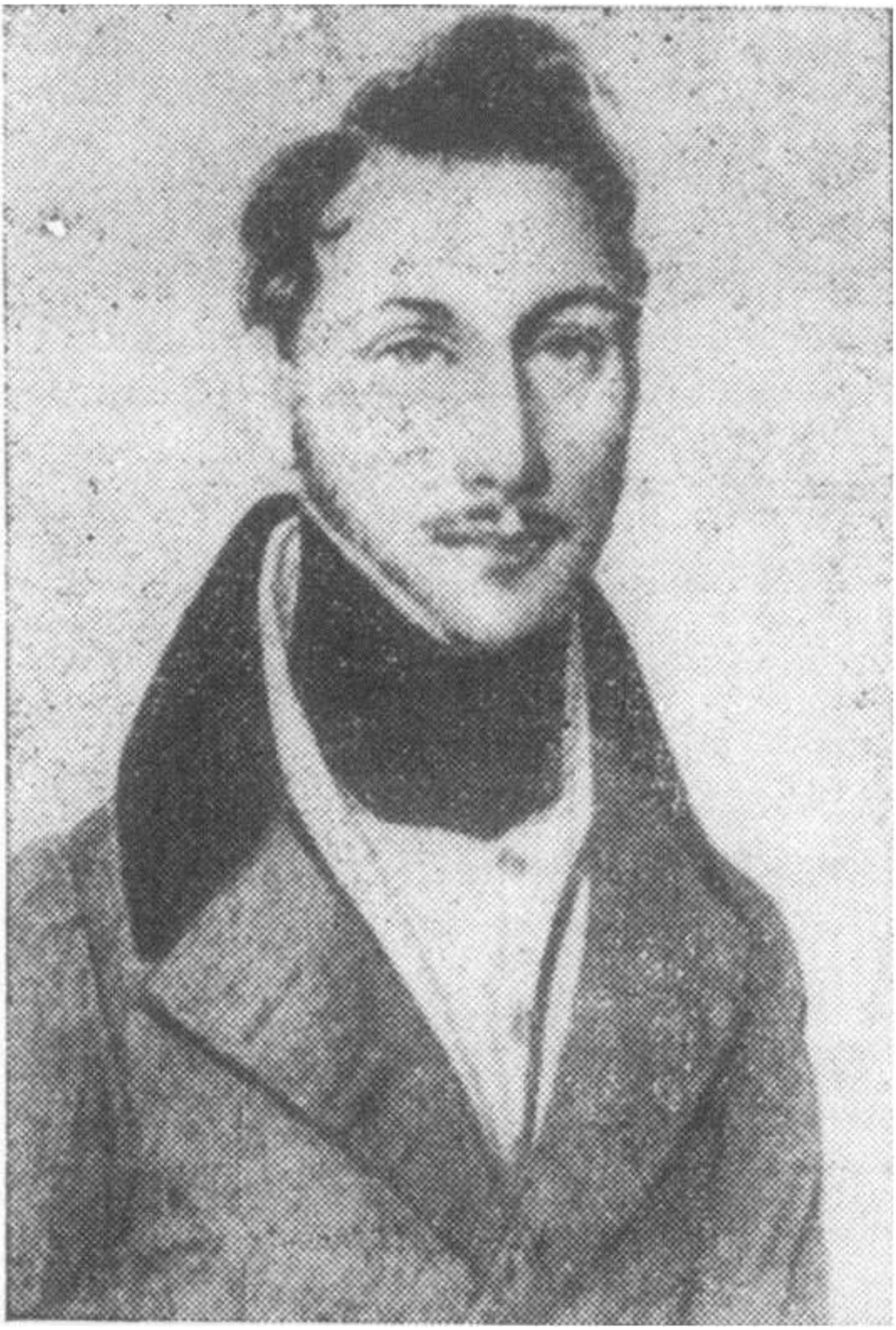

Picture 18. Dmitry Meyer.

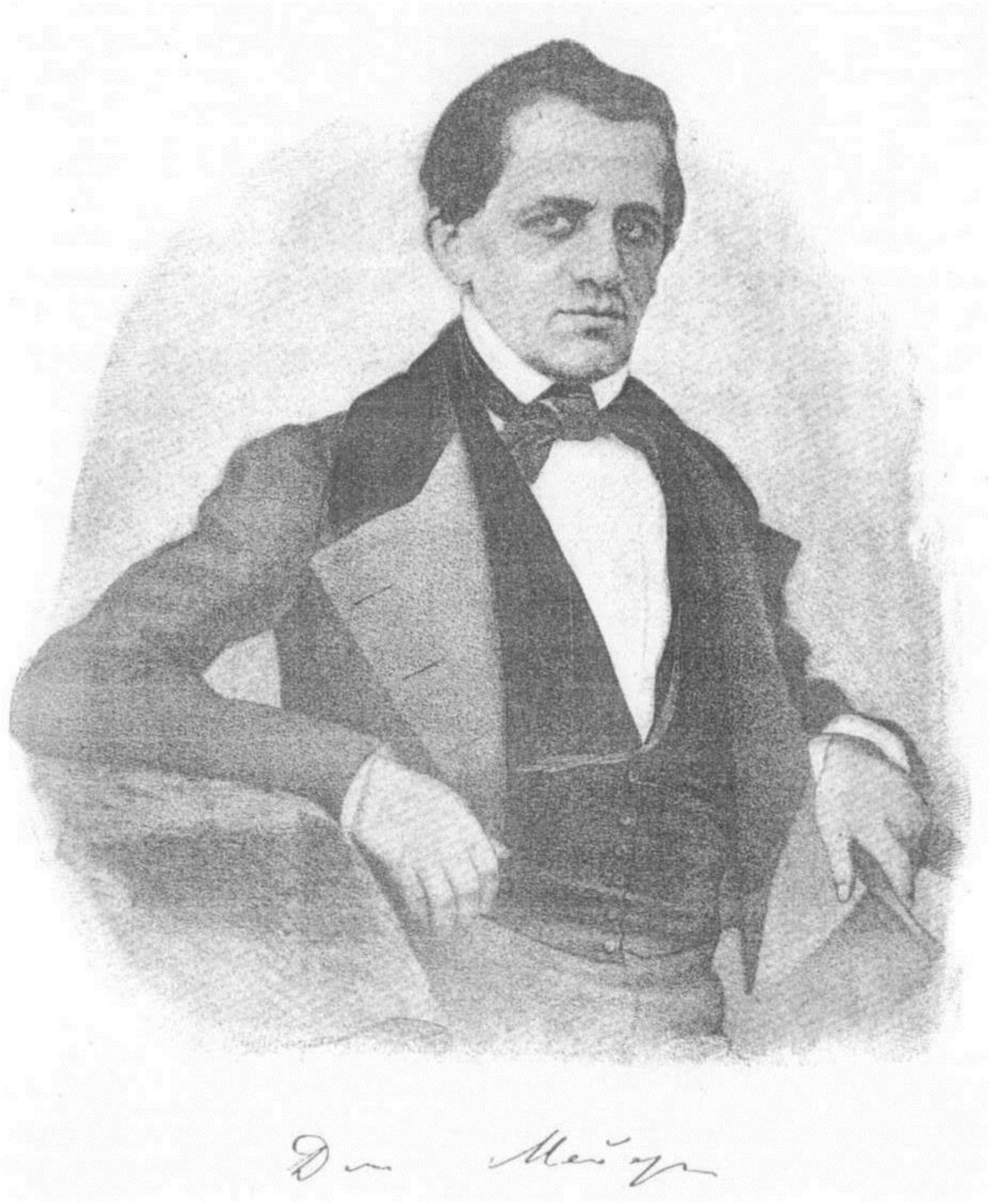